\documentclass[12pt]{amsart}
\usepackage{amsmath,amscd}
\usepackage{geometry,amsfonts,amssymb,amsthm,txfonts,pxfonts,amscd,url} 
\usepackage{algorithmicx,algpseudocode}
\usepackage{algorithm}
\usepackage{refcheck,graphicx,url}
\usepackage[section]{placeins}
\norefnames
\nocitenames
\def\struckint{\mathop{%
\def\mathpalette##1##2{\mathchoice{##1\displaystyle##2}%
  {##1\textstyle##2}{##1\scriptstyle##2}{##1\scriptscriptstyle##2}}%
\mathpalette
{\vbox\bgroup\baselineskip0pt\lineskiplimit-1000pt\lineskip-1000pt
\halign\bgroup\hfill$}
{##$\hfill\cr{\intop}\cr\diagup\cr\egroup\egroup}%
}\limits}
\newtheorem{theorem}{Theorem}[section]

\theoremstyle{remark}

\begin{document}
\title{Spectral experiments+}

\author{Igor Rivin}
\address{Department of Mathematics, Temple University, Philadelphia}
\email{rivin@temple.edu}
\date{\today}
\begin{abstract}
We describe extensive computational experiments on spectral properties of random objects - random cubic graphs, random planar triangulations, and Voronoi and Delaunay diagrams of random (uniformly distributed) point sets on the sphere). We look at bulk eigenvalue distribution, eigenvalue spacings, and locality properties of eigenvectors. We also look at the statistics of \emph{nodal domains} of eigenvectors on these graphs. In all cases we discover completely new (at least to this author) phenomena. The author has tried to refrain from making specific conjectures, inviting the reader, instead, to meditate on the data.
\end{abstract}
\thanks{The author would like to thank Peter Sarnak for interesting conversations, and the participants in the \url{http://mathematica.stackexchange.com} discussion group for help with matters computationa, and \url{http://mathoverflow.com} on help with pointers.
}
\maketitle
\tableofcontents
\section{Introduction}
There has been considerable attention devoted to spectral properties of random matrices, starting with Wishart's pioneering work in statistics (85 years ago) and Wigner's pioneering work in physics some seventy years ago. Since then, the subject has been central to mathematical physics, and following the work of Montgomery and Odlyzko (who discovered that spacings of zeros of the Riemann Zeta function on the critical line follow GUE statistics, to number theory, where the work of J. Keating's group (see, for example \cite{KeatingSurvey}) as well as N. Katz, P. Sarnak (see, for example, \cite{KatzSarnak}) and collaborators should be noted.

There has been considerable progress in the understanding of both global and local distribution of eigenvalues and singular values of symmetric, hermitian, and otherwise matrices, due to the efforts of L. Erdos, T. Tao, V. Vu, H-T. Yau, M. Rudelson, R. Vershynin, and many others. This author cannot hope to give a summary of literature in this note without slighting someone, so I will not even try. The interested reader is referred to the Oxford Handbook of Random Matrix Theory \cite{akemann2011oxford}(which, while quite new, is already out of date), or Mehta's classic \cite{mehta2004random}.

Over fifteen years ago, this author (jointly with D. Jakobson, S. Miller, and Z. Rudnick, see \cite{jakobson1999eigenvalue}) studied the distribution of eigenvalues of random \emph{regular} graphs. The bulk distribution of the eigenvalues of $k$-regular graphs has been understood since the work B. McKay \cite{McKay1981Eigen} - McKay showed that the probability density function of a large graph of degree $d$ approached
\[
f(x)=\begin{cases}
\frac{\sqrt{4(d-1)-x^2}}{2\pi(d^2-x^2)}, & \mbox{for $|x| \leq 2\sqrt{(d-1)}$}\\
0 & \mbox{otherwise.}
\end{cases}
\]
As $d$ tends to infinity, this approaches E. Wigner's semicircle law for the bulk distribution of the eigenvalues of random symmetric matrices.

McKay said nothing about the \emph{spacing distribution} of the eigenvalues, and this was the subject of our study. Our methods were quite primitive, but it was visually quite obvious that the spacing distribution "looked like" the GOE distribution, and there was visible level repulsion.

In the current study, we look at the following ensemble:
\begin{enumerate}
\item \label{goe} The base: GOE - random symmetric matrices with i.i.d. centered Gaussian entries.
\item Random cubic graphs.
\item Random planar 3-connected graphs (by E. Steinitz' theorem, these are the 1-skeleta of convex polyhedra.
\item Random Delaunay graphs, constructed as the $1$-skeleta of convex hulls of uniform samples of points on the sphere.
\item Random Voronoi triangulations (the planar duals of Delaunay graphs, as above
\end{enumerate}

For each of the ensembles we look at the bulk distribution of eigenvalues, the spacing distribution of adjacent eigenvalues, and also the \emph{localization} properties of the eigenvectors, as indicated by their $L^\infty$ norms, and the number of \emph{nodal domains} - for each eigenvector $v_\lambda,$ the nodal domains are the connected components of the subgraph induced by the set of vertices where $v_\lambda$ is positive (by symmetry, this is statistically the same as the subgraph induced by the vertices where $v_\lambda$ is negative.) All of these questions have been studied at great length in Case \ref{goe} - we certainly have fairly little to add, and that case is included mostly as the ``reference'' case.

\section{GOE}
\label{goesec}
We generated random symmetric matrices in the obvious way (we used matrices of size $3000 \times 3000,$ and in fact, in all cases these are the dimensions of our objects. Since the $L^\infty$ norm of the eigenvectors is very noisy, in order to generate somewhat informative graphs, we averaged the norm over intervals of length $0.02$ in eigenvalue space (all other localization graphs are averaged over the much smaller intervals of length $0.001,$ since the eigenvalues in the GOE model are much more dispersed.

\begin{figure}
\centering
\includegraphics[width=0.7\textwidth]{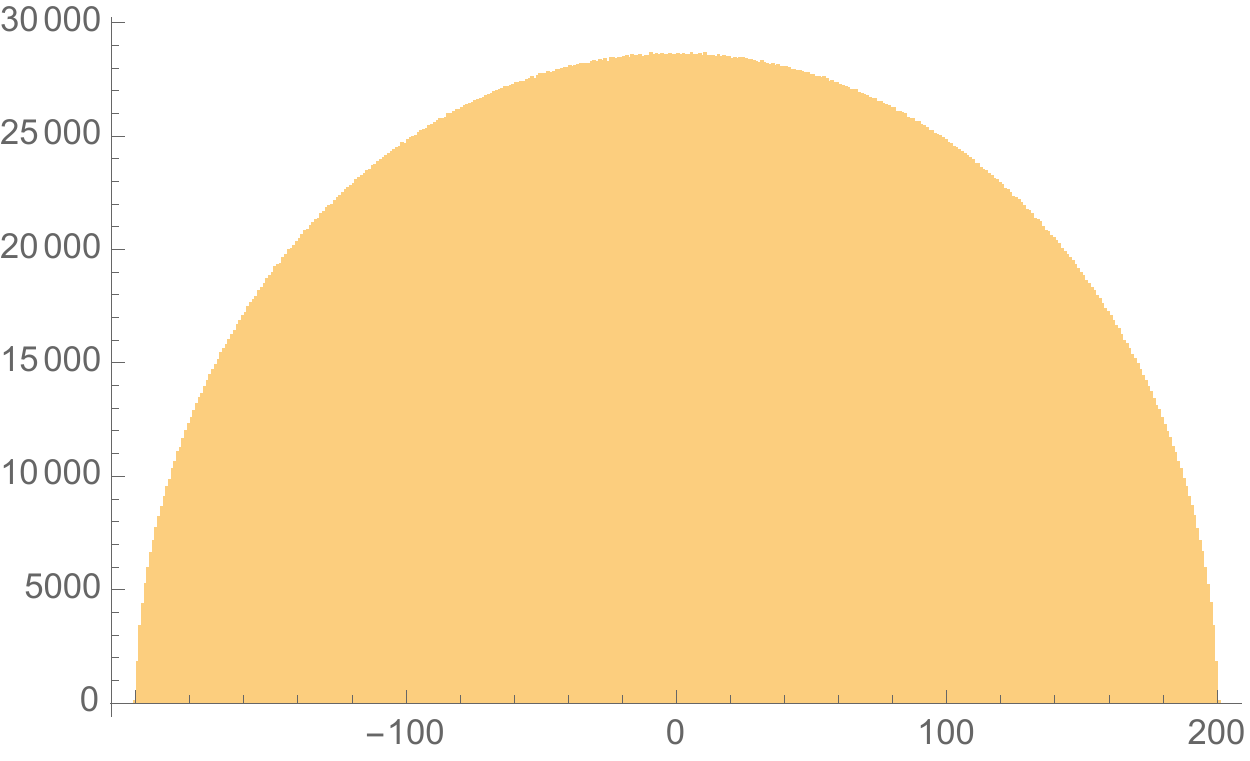}
\caption{\label{fig:wigner}Bulk Distribution of Eigenvalues of random symmetric matrices.}
\end{figure}

In any event, first, we have the bulk distribution in Figure \ref{fig:wigner}. To absolutely no-one's surprise, we see Wigner's semi-circle law in action.

\begin{figure}
\centering
\includegraphics[width=0.7\textwidth]{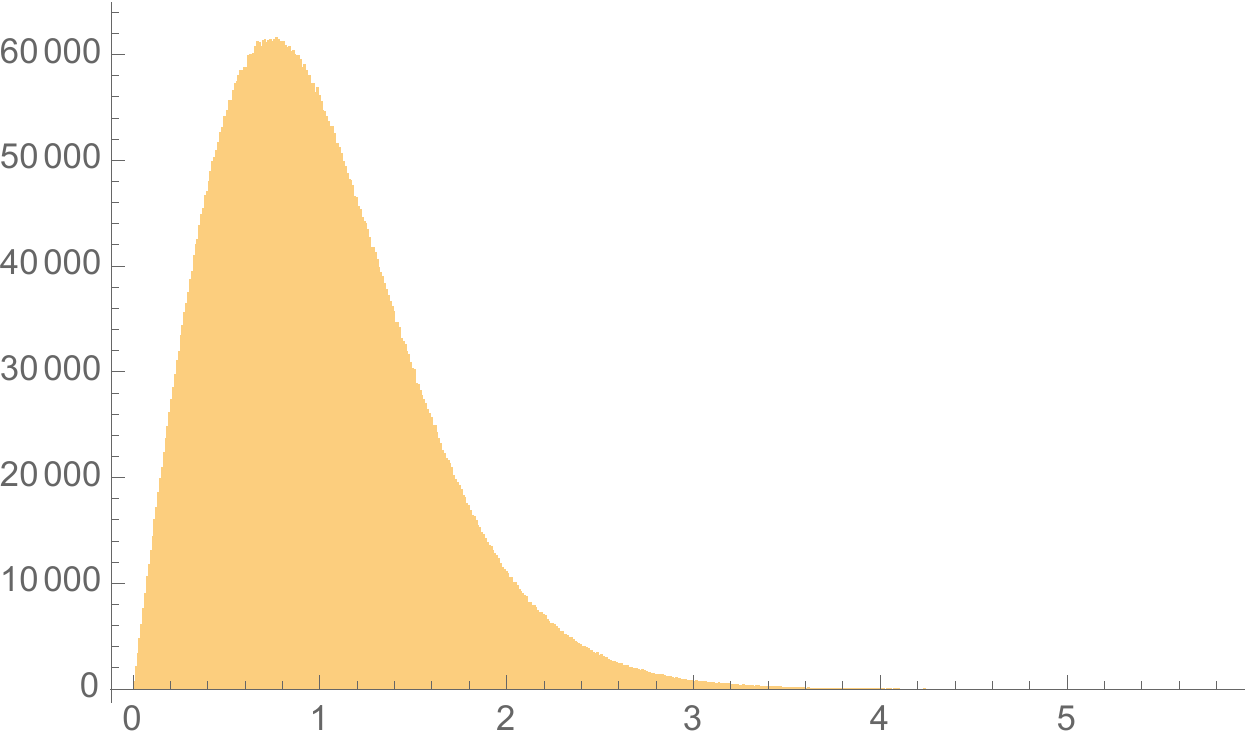}
\caption{\label{fig:spacings}Spacings Distribution of Eigenvalues of random symmetric matrices.}
\end{figure}

Next, we look at the spacing distribution (note: we do not do any "unfolding" by the Wigner density. We do, however, look only at spacings of the middle 80\% of the eigenvalues - we similarly throw out the edge eigenvalues for other spacing calculations.) We follow the physics custom of normalizing the mean spacing to $1.$ See Figure \ref{fig:spacings}.

\begin{figure}
\centering
\includegraphics[width=0.7\textwidth]{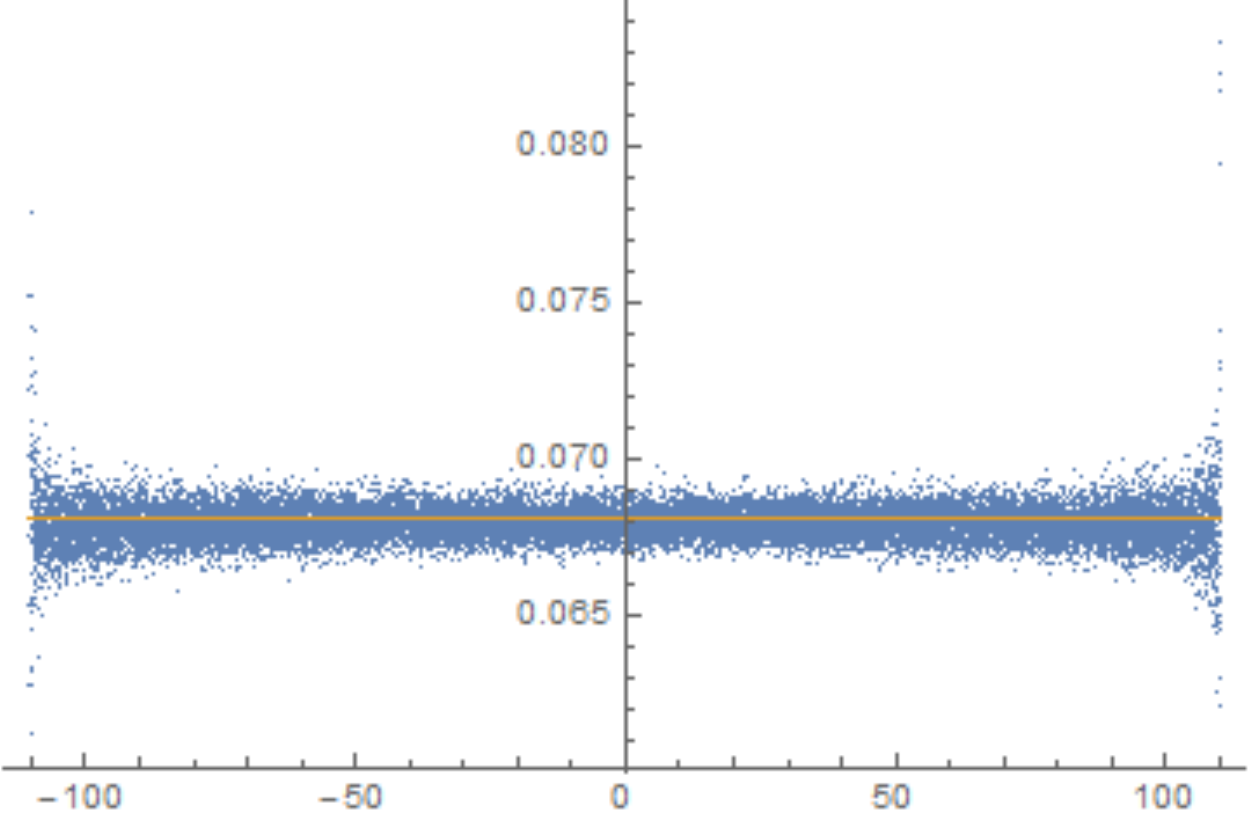}
\caption{\label{fig:localGOE}De-localization of eigenfunctions.}
\end{figure}

Lastly, we consider the localization properties of the eigenvectors of GOE matrices. In this case, from symmetry consideration (see, for example, the discussion in \cite{ErdosSchleinYau2009}), it seems intuitively clear that the $L^\infty$ norm of the eigenvector is roughly what it should be for random points on the $N$-dimensional sphere (where $N$ is the dimension of the matrix, and, indeed, this is borne out by experiment, see Figure \ref{fig:localGOE}.
Here (and in other similar graphs), the horizontal line represents the expected $L^\infty$ norm for a random point on the sphere. The only interesting aspect of this graph is the visible increase in variance toward the edge of the spectrum\footnote{in fact, in this, and other, localization graphs, the $L^\infty$ norms are average of an interval in eigenvalue space [here, the interval is of length $0.02,$ but in other localization graphs it is of length $0.001,$ since the eigenvalues are more closely spaced in the sparse cases]}
\section{Random cubic graphs}
\label{cubicgraphs}
\begin{figure}
\centering
\includegraphics[width=0.7\textwidth]{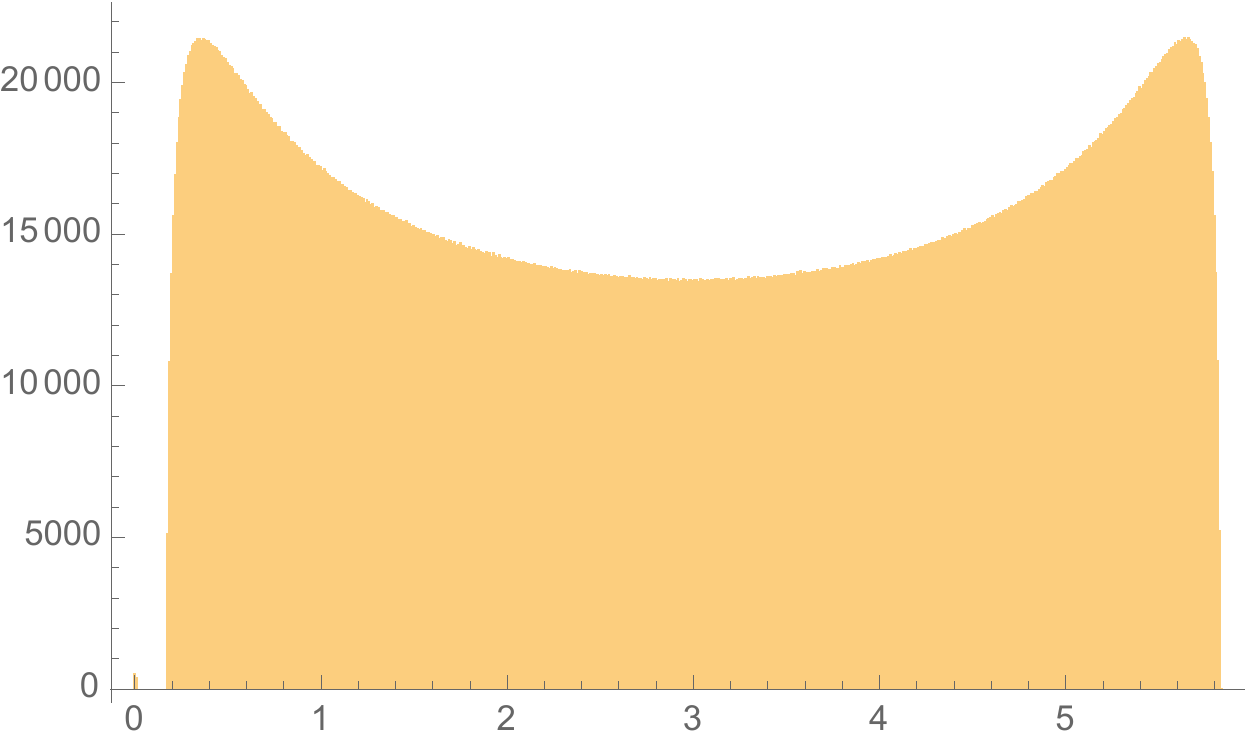}
\caption{\label{fig:cubicbulk}Bulk eigenvalue distribution for cubic graphs.}
\end{figure}

We revisit the random regular graphs of \cite{jakobson1999eigenvalue}, generated using B. Bollobas' algorithm. Not surprisingly, the bulk distribution of eigenvalues is in perfect agreement with McKay's law, see Figure \ref{fig:cubicbulk}.

\subsection{Spacings}
\begin{figure}
\centering
\includegraphics[width=0.7\textwidth]{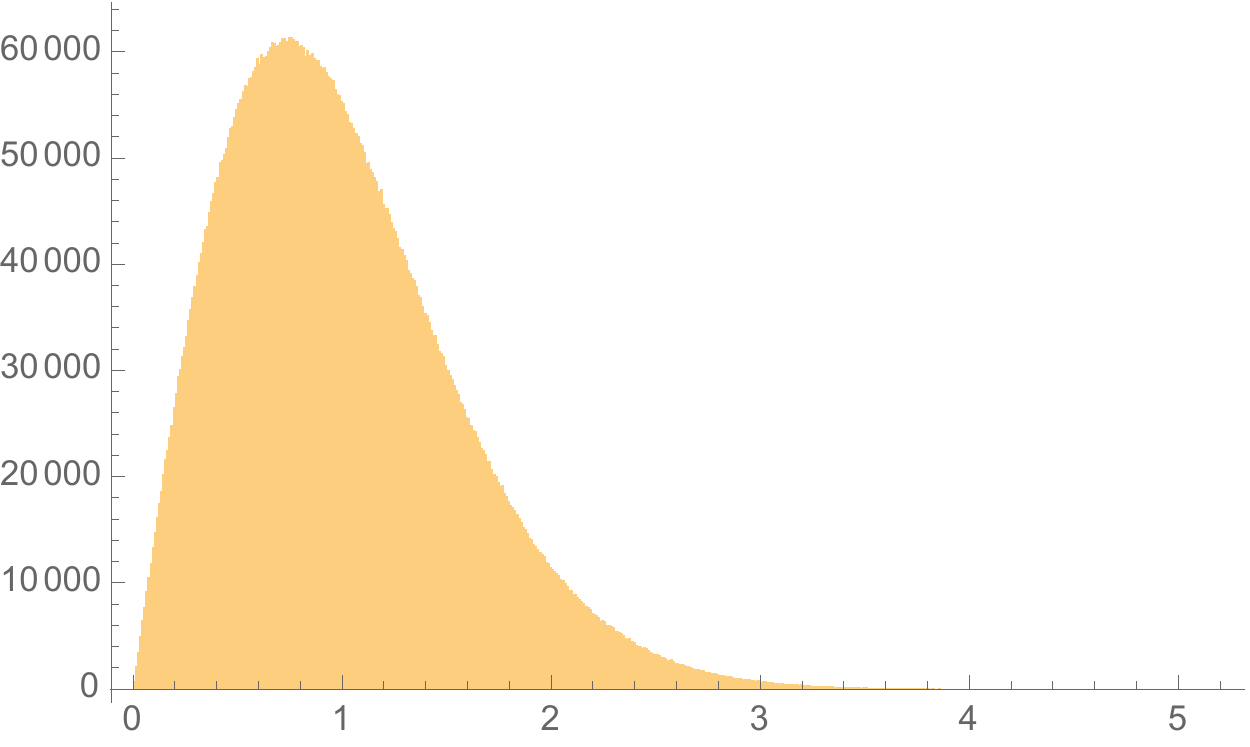}
\caption{\label{fig:cubicspacings}Spacing distribution for cubic graphs.}
\end{figure}

Next, we study eigenvalue spacings, as before plucked from the middle 80\% of the eigenvalues, and as before unfolded, in Fig.~\ref{fig:cubicspacings}

Now, the spacing distribution looks awfully like the GOE spacings figure \ref{fig:localGOE}. Kolmogorov-Smirnov's test rejects this out of hand, but a quantile plot \ref{fig:cubicvsGOEQQ} shows us what is happening:
\begin{figure}
\centering
\includegraphics[width=0.7\textwidth]{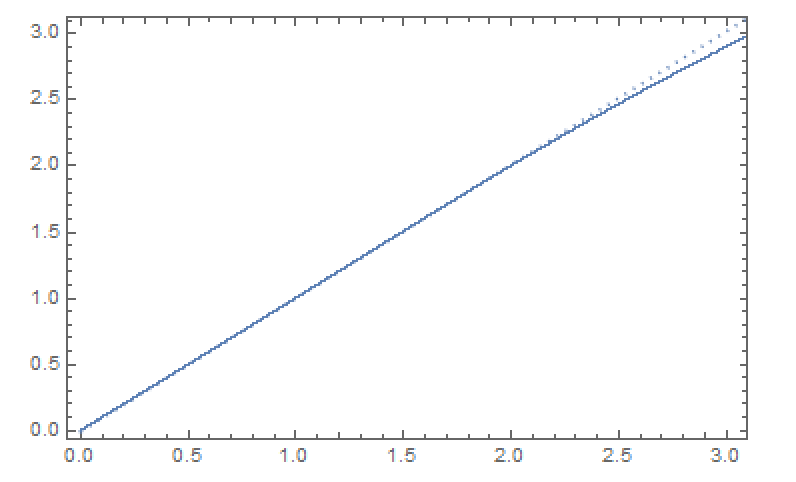}
\caption{\label{fig:cubicvsGOEQQ}Quantile Plot of cubic spacing vs GOE spacing.}
\end{figure}
The reader will see that the there is a small, but quite noticeable divergence at the high end of the distribution, so the tail for graphs is a little less fat. The excellent agreement through most of the range is still consistent with the hypothesis that the limiting distributions are the same.

\subsection{Localization}
\begin{figure}
\centering
\includegraphics[width=0.7\textwidth]{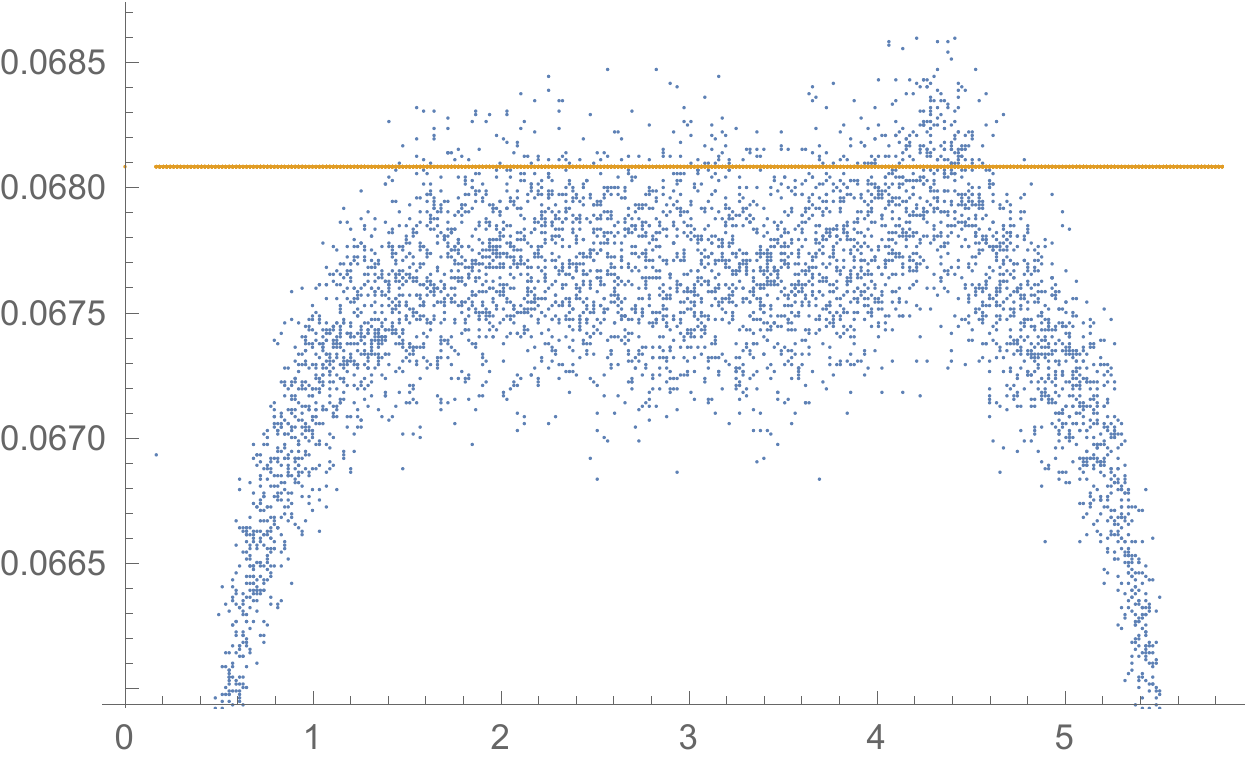}
\caption{\label{fig:cubiclocal}Spacing distribution for cubic graphs.}
\end{figure}

Finally, we look at the $L^\infty$ norms of the eigenvectors of random cubic graphs, and here we are in for a bit of a surprise -- see Figure \ref{fig:cubiclocal}. 
\begin{itemize}
\item 
First, we see that there is \emph{almost} (but not quite) a symmetry around the midpoint of the spectrum (if we had looked at the adjacency, instead of the Laplacian, matrix) the symmetry would be more apparent, but there is a visible spike around $3 + \sqrt{2}.$
\item
Secondly, there is a very visible collapse in the $L^\infty$ near the edges of the spectrum, and there appears to be a phase transition around $3 \pm \sqrt{2}.$ 
\item
Thirdly (and most surprisingly) it seems that with the exception of the spike, the eigenfunctions on cubic graphs are \emph{more uniformly} distributed that could be expected from random functions.
\end{itemize}
Of course, on the level of eigenfunctions, random cubic graphs look nothing at all like GOE matrices, and this author has no explanation for the behaviors observed.
\section{Random planar cubic three-connected graphs}
\label{steinitzsec}
It is a well-known theorem of E.~Steinitz that a graph is the $1$-skeleton of a convex polyhedron in $\mathbb{E}^3$ if and only if it is planar and $3$-connected (for a nice proof, see G. Ziegler's book \cite{ziegler1995lectures}). Luckily for us, the \emph{Boltzmann sampler} of P.~Duchon, P.~Flajolet, G.~Louchard, and G.~Schaeffer \cite{duchon2004boltzmann} can be used to generate random such graphs uniformly.footnote{in order for this to be efficient, the graphs need to be close to [but not exactly equal to] the desired size, however, this is easy to account for in the statistics} This has been implemented by Schaeffer, and distributed as \texttt{PlanarMap}. It should be noted that the results in this section are clearly related to the geometry of the \emph{Brownian Map} - see \cite{gall2014random}, and references therein.
So, with this in hand, we can see what we get:
\subsection{Bulk distribution of eigenvalues}
\begin{figure}
\centering
\includegraphics[width=0.7\textwidth]{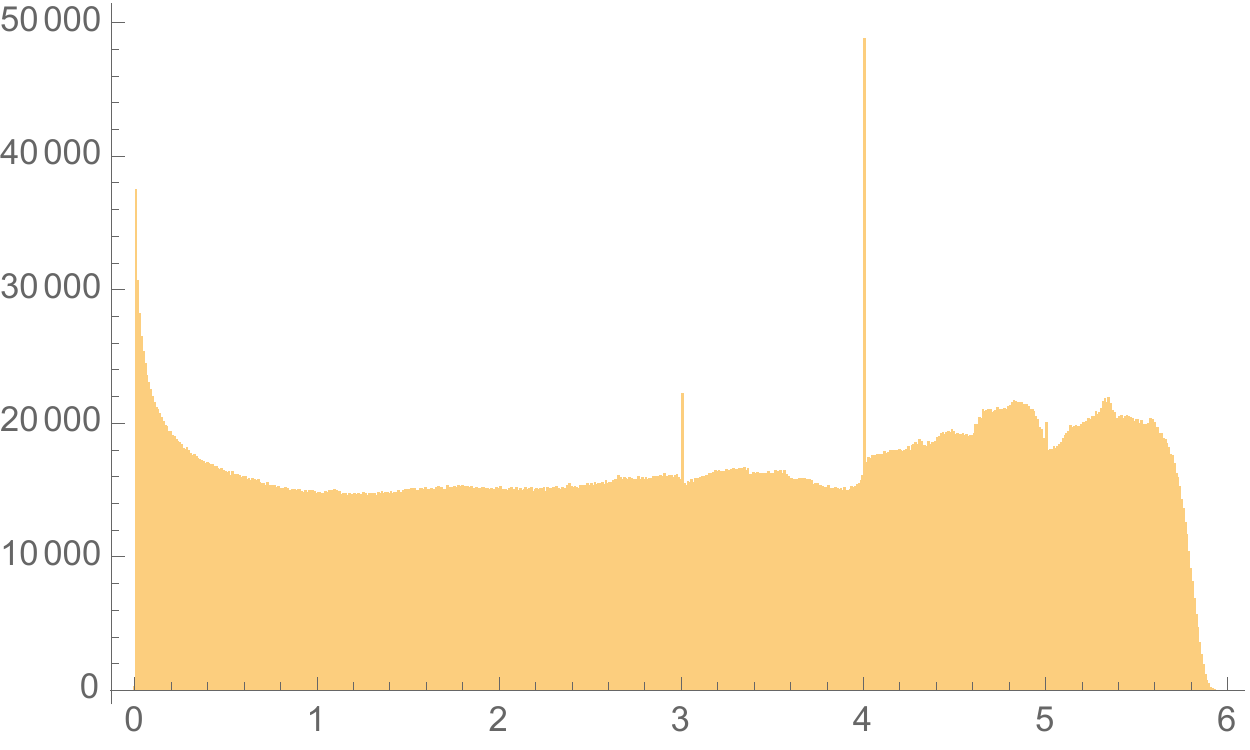}
\caption{\label{fig:randomPGBulk}Bulk distribution of random cubic map eigenvalues}
\end{figure}

First, the bulk distribution of the eigenvalues -- see Figure \ref{fig:randomPGBulk}.
We note a number of interesting features:
\begin{itemize}
\item there is a noticeable spike in density at $0$ -- the fact that the density at $0$ is nonzero is quite well-known (planar graphs have no spectral gap), but the spike is  a little  surprising.
\item Peculiarly, there are sharp spikes at $3, 4, 5,$ which this author  is at a loss to explain.
\item Between $1$ and $3$ the distribution appears to be close to uniform.
\item The ``high energy'' levels seem to be relatively irregular.
\end{itemize}
\subsection{Spacings}
\begin{figure}
\centering
\includegraphics[width=0.7\textwidth]{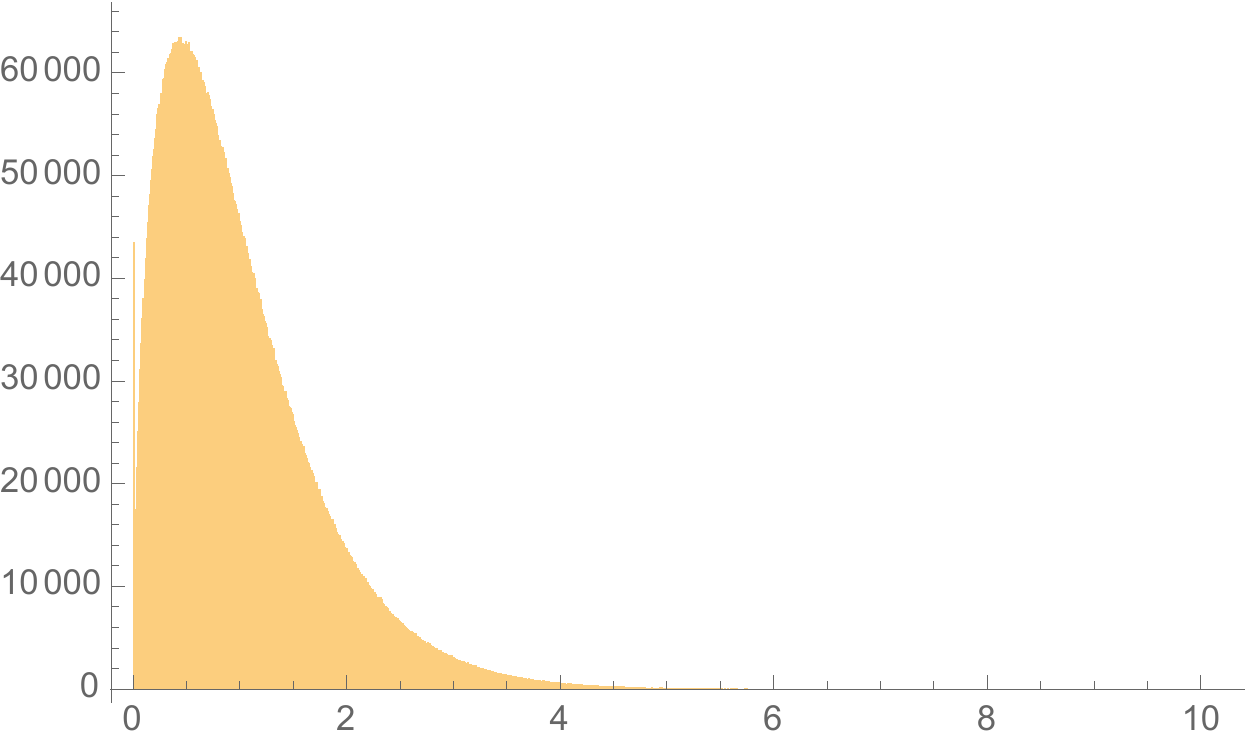}
\caption{\label{fig:randomPGspacings}Spacing distribution for random $3$-connected cubic planar maps}
\end{figure}

\begin{figure}
\centering
\includegraphics[width=0.7\textwidth]{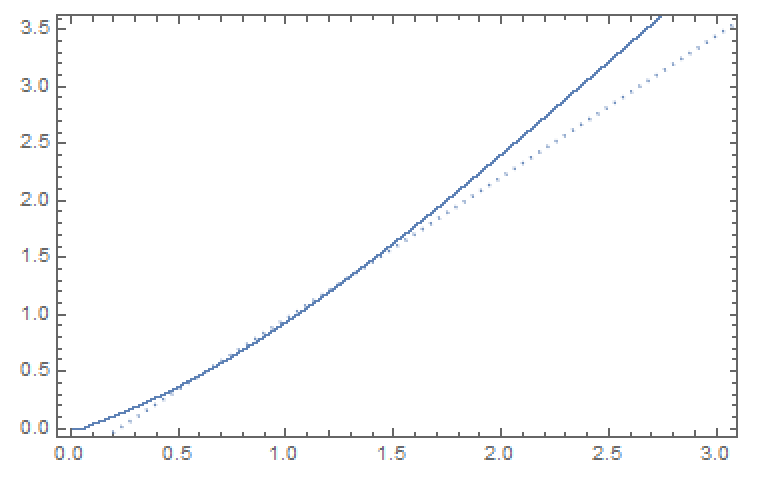}
\caption{\label{fig:randomPGGOEQQ}Quantile plot of planar cubic spacings vs GOE}
\end{figure}

Next, the spacings. First, the graph - see Figure \ref{fig:randomPGspacings}.
A quick look at the figure reveals that:
\begin{itemize}
\item There is still level repulsion.
\item The tails are much fatter than in the GOE distribution.
\end{itemize}
We can quantify these observations with a quantile plot, as before - see Figure \ref{fig:randomPGGOEQQ}.
\subsection{Localization}
\begin{figure}
\centering
\includegraphics[width=0.7\textwidth]{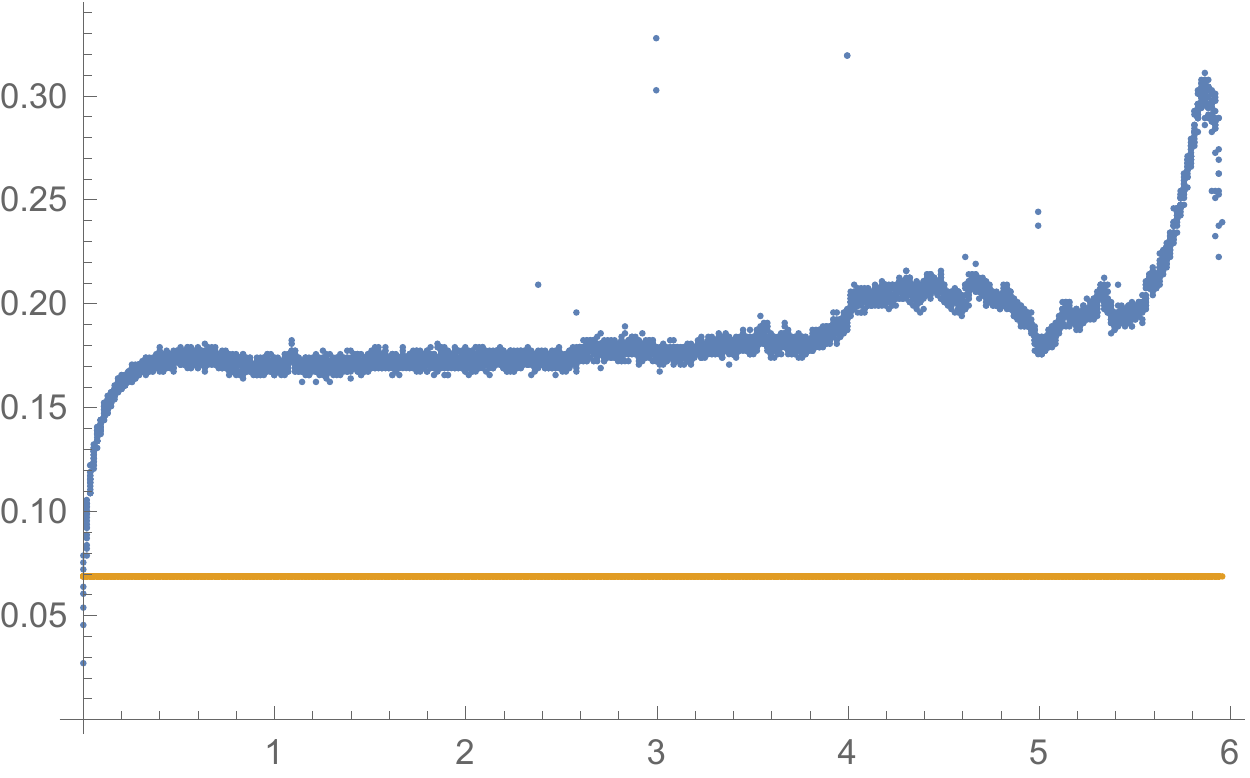}
\caption{\label{fig:planarlocal}$L^\infty$ norm of eigenvectors for random cubic $3$-connected planar maps}
\end{figure}

The graph of the $L^\infty$ norm of the eigenvectors for planar cubic maps (vs what it would be if the eigenvectors were behaving like random functions) can be found in Figure \ref{fig:planarlocal}.
Various observations can be made:
\begin{enumerate}
\item The $L^\infty$ norms are clearly consistent with a lot of structure.
\item There are clear phases, which map onto the phases visible  in the bulk eigenvalue distribution:
\item There is clearly something special going on at the magic integer values $3, 4, 5.$
\item The graph is flat in the same region between $1$ and $3$ where the eigenvalue density appears flat.
\item There is a noticeable spike in the $L^\infty$ norm around what looks like  $3 + 2\sqrt{2}$ in eigenvalue space.
\end{enumerate}
\section{Delaunay and Voronoi triangulations of planar point sets}
Given a collection of points $\mathcal{P}$ on the unit sphere $\mathbb{S}^2 \subset \mathbb{E}^3,$ we can form their convex hull $C(\mathcal{P}).$ If we pick some $p\in \mathcal{P}$ as the north pole, stereographic projection from $p$ transforms $C(\mathcal{P})$ into the \emph{Delaunay triangulation} of the image point set (for more on this construction, see \cite{rivin1994euclidean}) - of course, the term ``triangulation'' is not strictly speaking justified, since it is possible that some faces of the convex hull are not triangles, but this is a probability $0$ event, which never happens when the points are chosen at random (as they are in these experiments). The combinatorics of Delaunay triangulations is restricted - the precise nature of the restriction was determined in \cite{rivin1996characterization} and, in higher genus,in \cite{rivin2003combinatorial}.

The dual of a Delaunay triangulation is a \emph{Voronoi diagram} (see H. Edelbrunner's beautiful book \cite{edelsbrunner2001geometry} for all you ever wanted to know on the subject). The Voronoi diagram (being the dual of a triangulation) is a cubic graph, so our results on Voronoi diagrams are more compatible with what came before.

\section{Voronoi diagrams}
\label{voronoisec}
\begin{figure}
\centering
\includegraphics[width=0.7\textwidth]{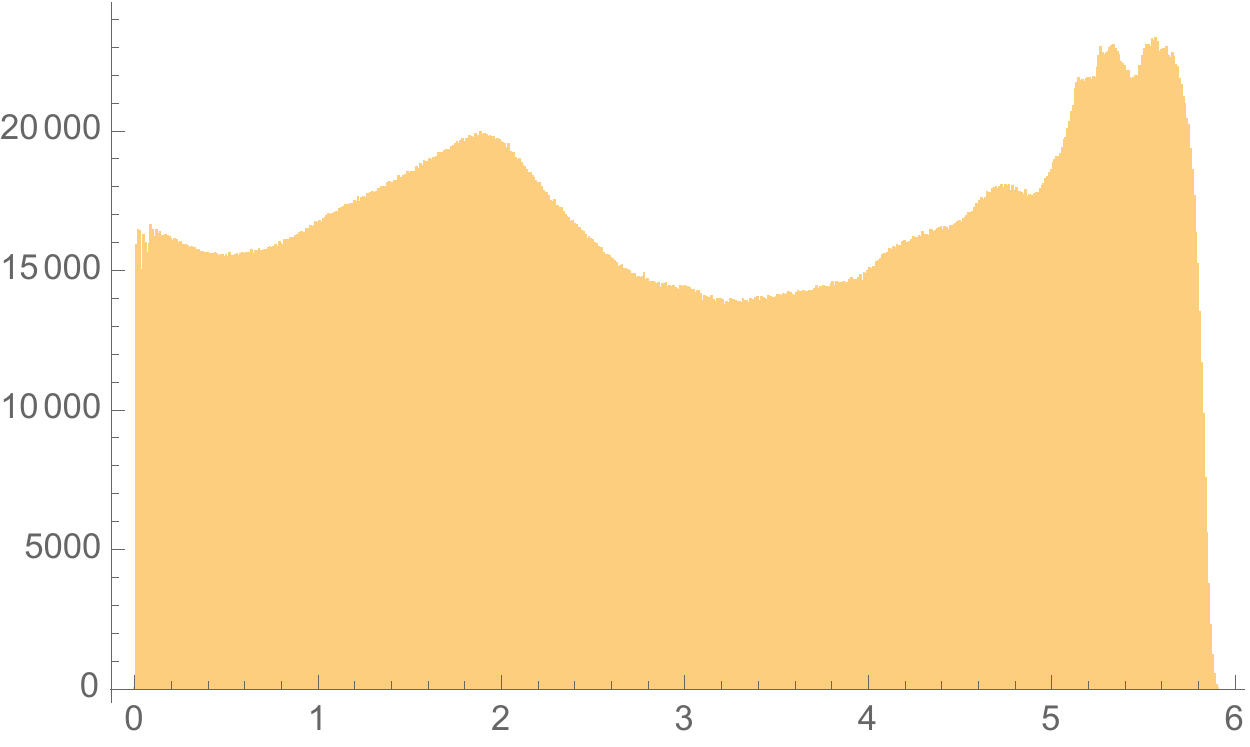}
\caption{\label{fig:voronoiBulk}Bulk distribution of eigenvalues for random Voronoi diagrams}
\end{figure}

First, the bulk distribution -- see Figure \ref{fig:voronoiBulk}.
The bulk distribution looks as bit messy, but many phases can be distinguished. 
\begin{itemize}
\item{Phase 1} The low-energy phase (with the peak at 0).
\item{Phase 2} The mid-low-energy phase, between roughly $0.7$ and $1.8.$ Here, the spectral density grows linearly, which would seem to correspond to the range in which the combinatorial Lalpacian approximates the geometric Laplacian (so the linear growth is what is predicted by Weyl's law.
\item{Phase 3} Mid-high-energy phase (until around $4.7$.
\item{Phase 4} High energy phase from $4.7$ to (almost) $6.$
\end{itemize}

With the exception of the mid-low-energy phase, this distribution seems somewhat mysterious at the moment.

\subsection{Spacings}
\begin{figure}
\centering
\includegraphics[width=0.7\textwidth]{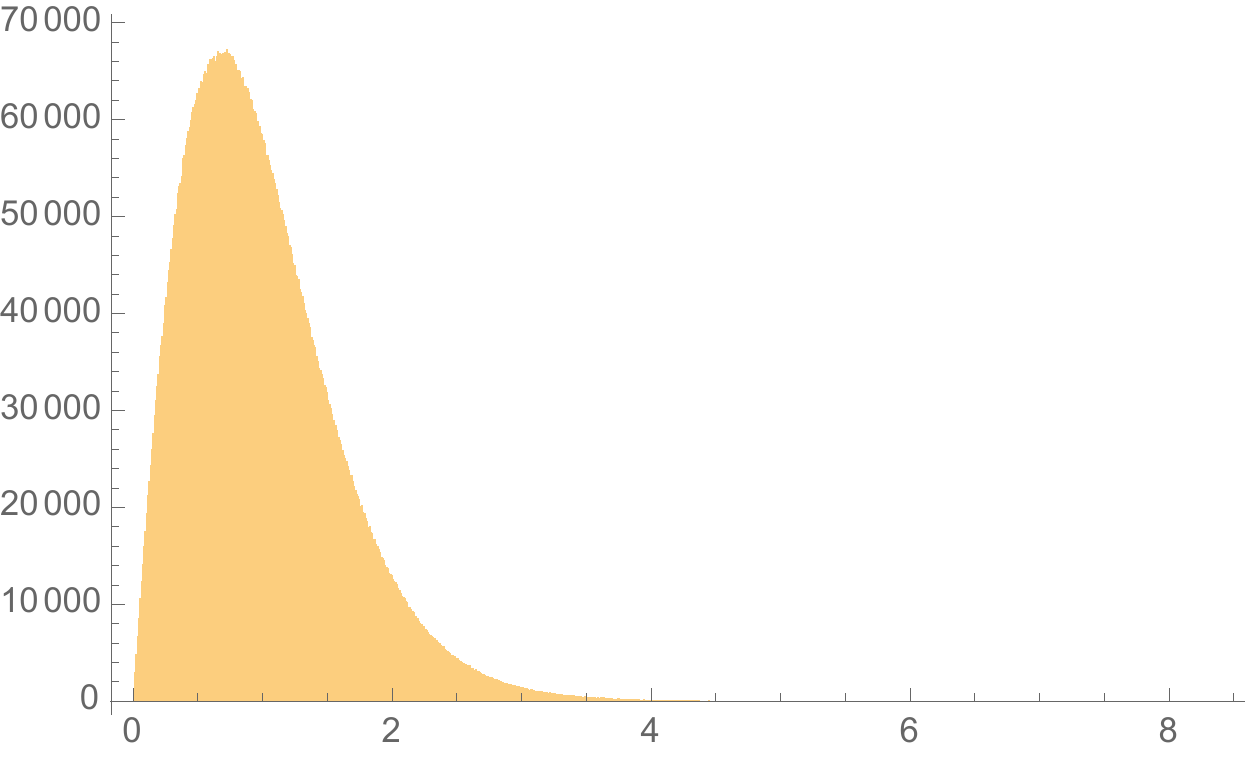}
\caption{\label{fig:voronoiSpacings}Bulk distribution of eigenvalues for random Voronoi diagrams}
\end{figure}

\begin{figure}
\centering
\includegraphics[width=0.7\textwidth]{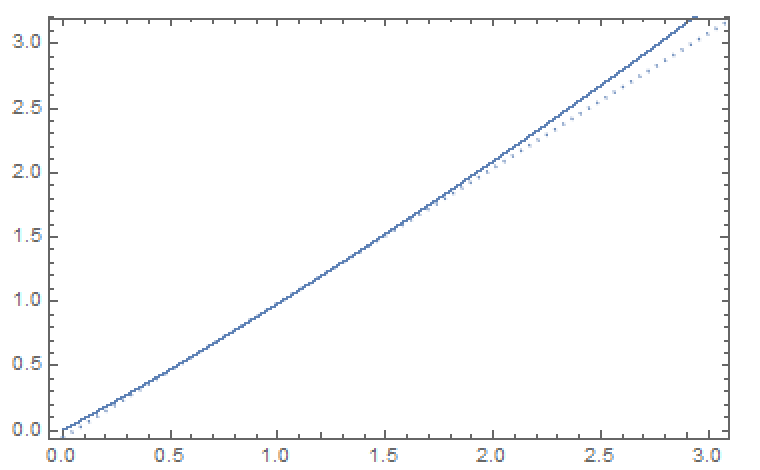}
\caption{\label{fig:voronoiGOEQ}Bulk distribution of eigenvalues for random Voronoi diagrams}
\end{figure}

The spacing graphs for random Voronoi diagrams (Figure \ref{fig:voronoiSpacings} looks like every other spacing graph we have seen - there is clearly level repulsion, but the distribution is \emph{not} that of GOE spacings, as can be seen in the quantile plot (note, however, that the distribution is much closer to GOE than the random planar cubic maps distribution) - Figure \ref{fig:voronoiGOEQ}. In fact, it may be reasonable to conjecture that the limiting spacing is GOE.
\subsection{Eigenfunction localization}
The localization picture for random Voronoi diagrams bears some similarity to the picture for random planar maps - see Figure \ref{fig:voronoiLocal}:
\begin{figure}
\centering
\includegraphics[width=0.7\textwidth]{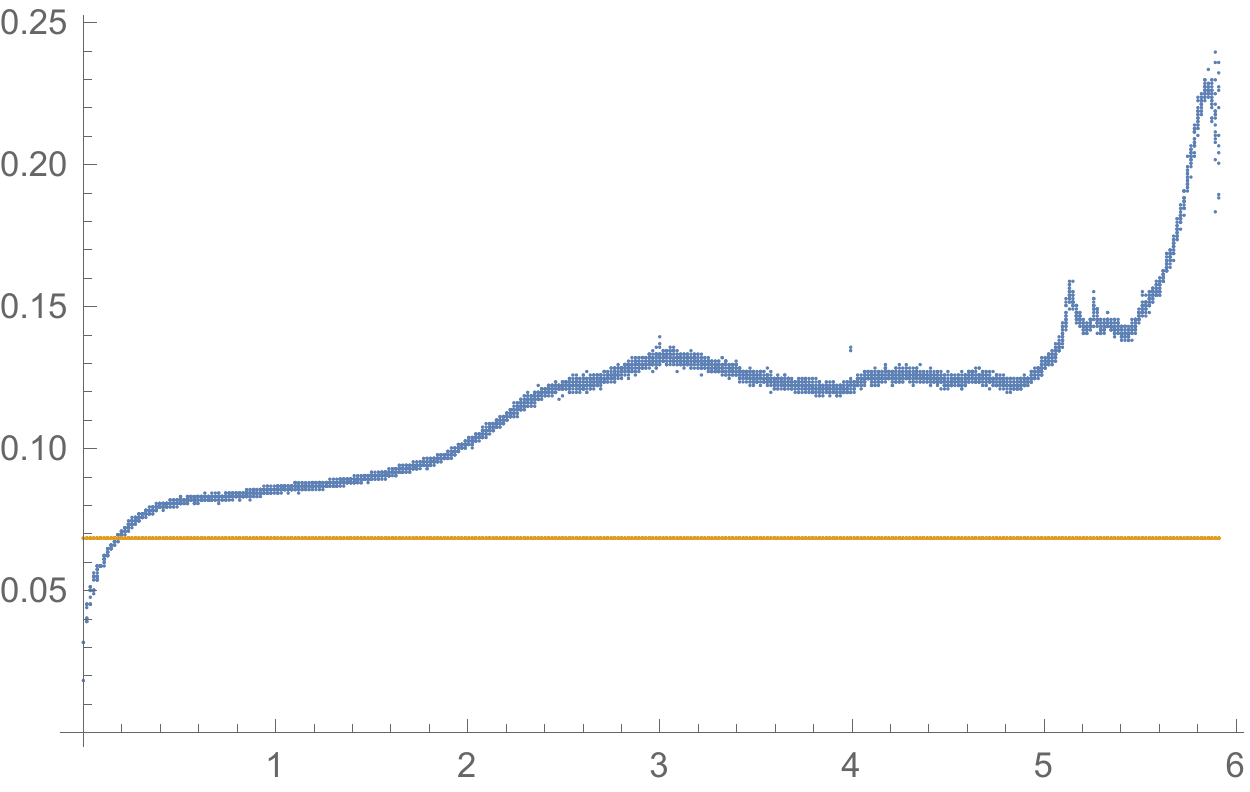}
\caption{\label{fig:voronoiLocal}$L^\infty$ norms of eigenfunctions for random Voronoi diagrams}
\end{figure}
\begin{enumerate}
\item The low energy eigenfunctions are very uniformly distributed.
\item the mid-energy eigenfunctions are clearly more localized than what would be predicted by randomness, but their $L^\infty$ norms seem quite stable.
\item There is a high energy spike (and subsequent decline) around $5.8$ in eigenvalue coordinates.
\end{enumerate}
\section{Delaunay triangulations}
\label{delaunaysec}
Delaunay triangulations (more precisely, their 1-skeleta) are not cubic graphs, so we expect their statistics to be different from our other graph ensembles. However, they \emph{are} dual to Voronoi diagrams, so we would expect to see considerable similarities between the two ensembles.
\subsection{Bulk distribution}
\begin{figure}
\centering
\includegraphics[width=0.7\textwidth]{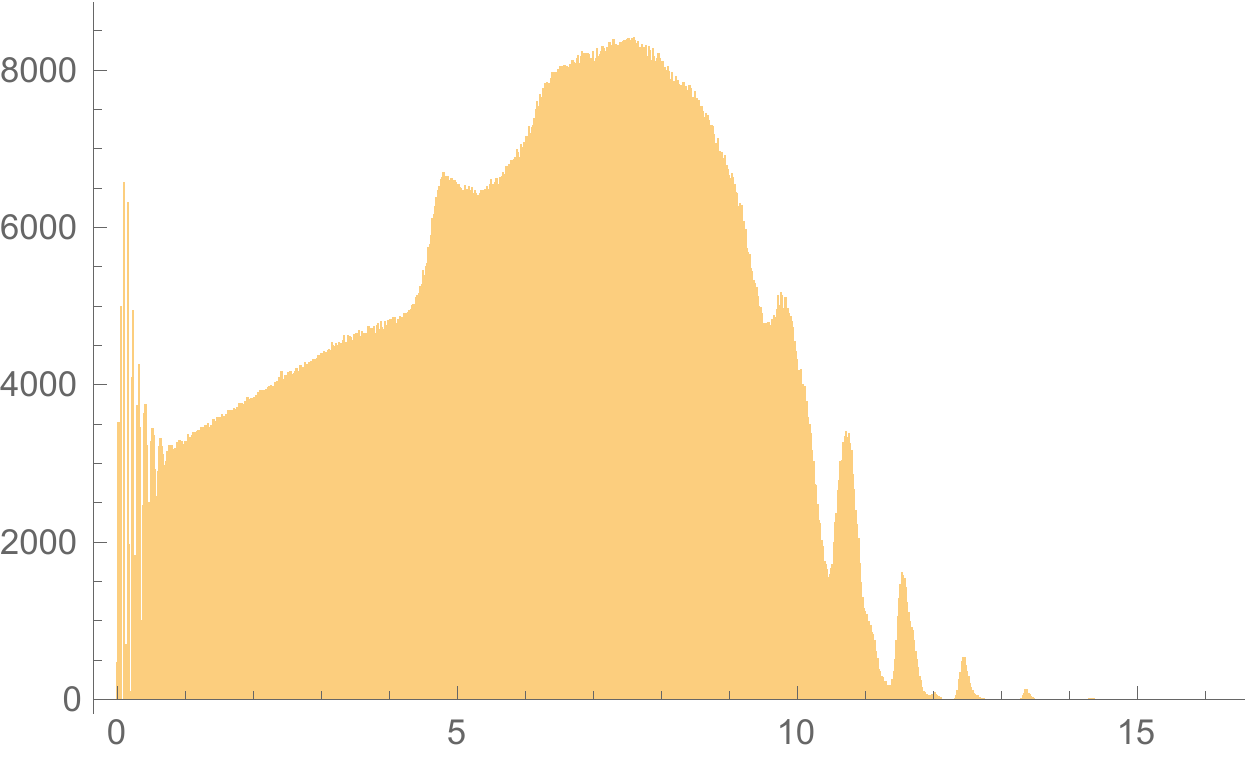}
\caption{\label{fig:delaunayBulk}Bulk spectral distribution for Delaunay triangulations of random point sets}
\end{figure}

The bulk distribution (Figure \ref{fig:delaunayBulk}) does present \emph{some} similarities to the Voronoi picture, together with major differences.
\begin{itemize}
\item A major similarity is the linear growth of spectral density int the mid-low energy regime (between abscissa $1$ and $5,$ approximately).
\item A minor similarity is the high density in very low energy regime.
\item Another is the continuing increase in spectral density through the mid-high energy.
\item A major \emph{difference} is the very noisy behavior at both very low and very high energies, with massive spikes in both regimes, which seem to indicate some sort of "phase locking" (that is, just as for random planar maps, the integral values $3, 4, 5$ seem special, here it seems that there is some sort of discreteness taking place.
\end{itemize}
\subsection{Spacing distribution}
\begin{figure}
\centering
\includegraphics[width=0.7\textwidth]{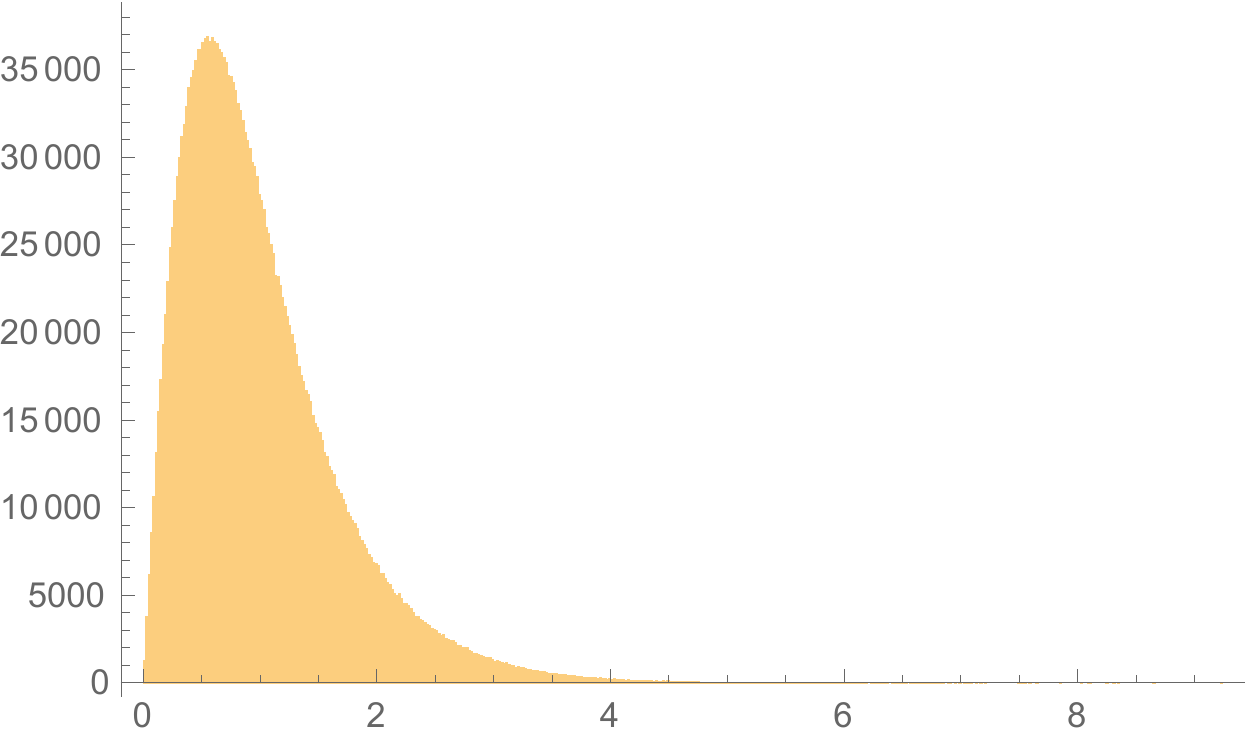}
\caption{\label{fig:delaunaySpacings}Spectral spacing distribution for Delaunay triangulations of random point sets}
\end{figure}

\begin{figure}
\centering
\includegraphics[width=0.7\textwidth]{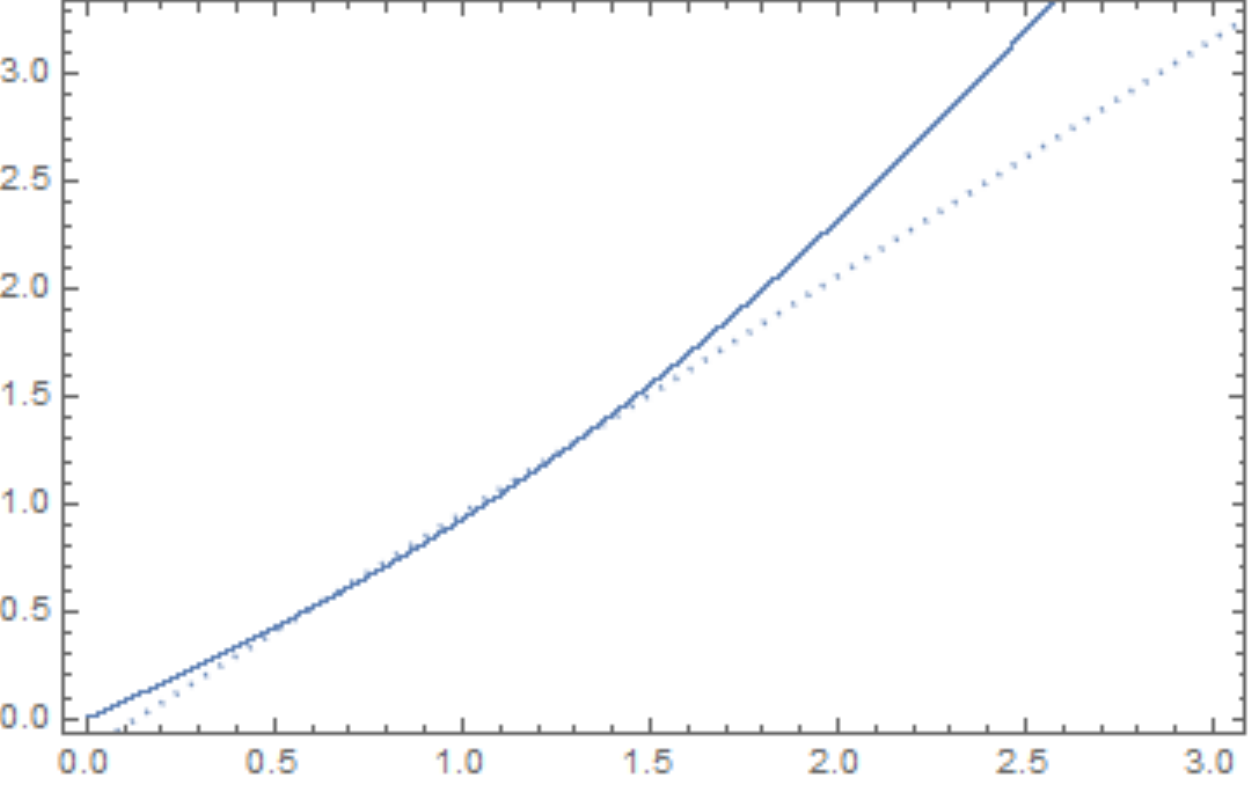}
\caption{\label{fig:delaunayQQ}Random Delaunay spacings vs GOE spacings quantile plot}
\end{figure}

The spacing distribution looks the same as ever (Figure \ref{fig:delaunaySpacings}). It visually appears to have fatter tails than GOE, and the quantile plot (Figure \ref{fig:delaunayQQ}) bears this out - indeed this is the worst fit yet to the GOE spacings distribution, though again, it seems hard to deny the presence of level repulsion.
\subsection{Locality of eigenfunctions}
\begin{figure}
\centering
\includegraphics[width=0.7\textwidth]{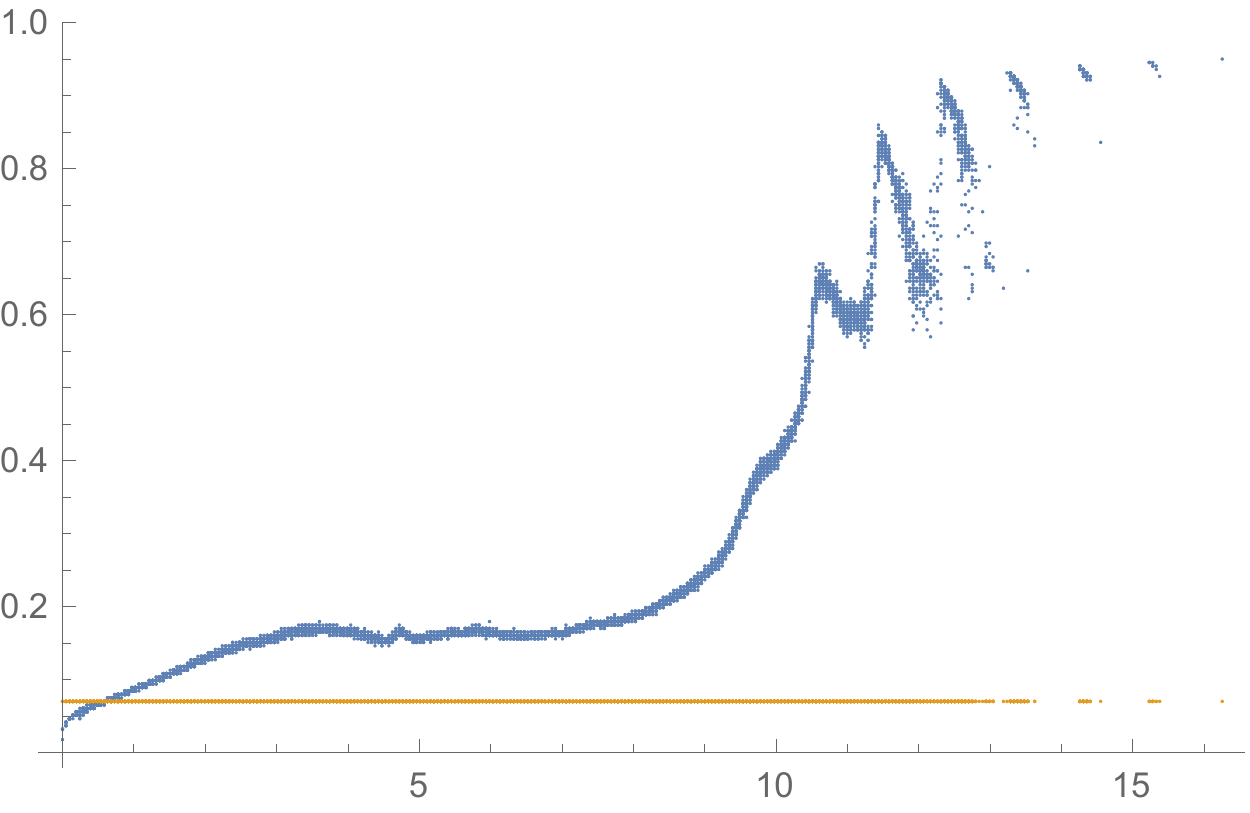}
\caption{\label{fig:delaunayLocal}Localization of eigenfunctions for Delaunay triangulations of random point sets}
\end{figure}

Finally, we look at the localization of the eigenvectors in the random Delaunay ensemble, see Figure \ref{fig:delaunayLocal}.
Here we see a major surprise:
\begin{itemize}
\item In the low-energy regime through the mid-low energy, we have a picture similar to that in Figure \ref{fig:voronoiLocal}, and in fact, the $L^\infty$ norm stays close to constant from around $3$ to $5.$ 
\item but then it spikes sharply, and goes to full-on localization (the $L^\infty$ norm tops out at $1$).
\end{itemize}
One might be tempted to explain the high-energy localization by the low correlation length (the graphs are planar) which makes the graphs look "almost banded", but this reasoning is clearly invalid, since if it worked, we would have the same results in the random planar map case and in the Voronoi diagrams case (especially the latter, since topologically there should not be that much difference between a cell complex and its planar dual), but they clearly do not.
\section{Nodal domains}
In this section we collect our results on nodal domains of eigenvectors in the various models of graphs we looked at, and a few more. The reader should note that the figures below were obtained by looking at a \emph{single random graph} of one of our ensembles, which indicates that the nodal domain statistics exhibit strong concentration properties. It should also be noted that one model we do \emph{not} look at is that of the Erdos-Renyi ranodm graphs $G(n, p).$ This has been studied by Dekel, Lee, and Linial \cite{dekel2011eigenvectors}, where they showed tha there in this model there are always $O_p(1)$ nodal domains.

For random regular graphs, there are results which have the flavor of Courant's theorem, to the effect that the number of nodal domains grows at most linearly with the ordinal number of the eigenvalue. 

First, let us recall the theorem of Lin, Lippner, Mangoubi, Yau \cite{lin2013nodal}:
\begin{theorem}[\cite{lin2013nodal}]
\label{yauthm}
Let the eigenvalues of $G$ be ordered in inreasing order, and let the maximal degree of $G$ be bounded above by $d.$ Then, the number $N(\lambda_n)$ of Nodal domains for the $n$-th eigenvalue is at most
\[
N(\lambda_n) \leq (d-1) n
\]
\end{theorem}
Since the number of nodal domains cannot be bigger than the number of vertices, if $G$ is actually $d$-regular, then Theorem \ref{yauthm} is only non-vacuous for $n<V(G)/(d-1).$
\subsection{Random regular graphs}
\begin{figure}
\centering
\includegraphics[width=0.7\textwidth]{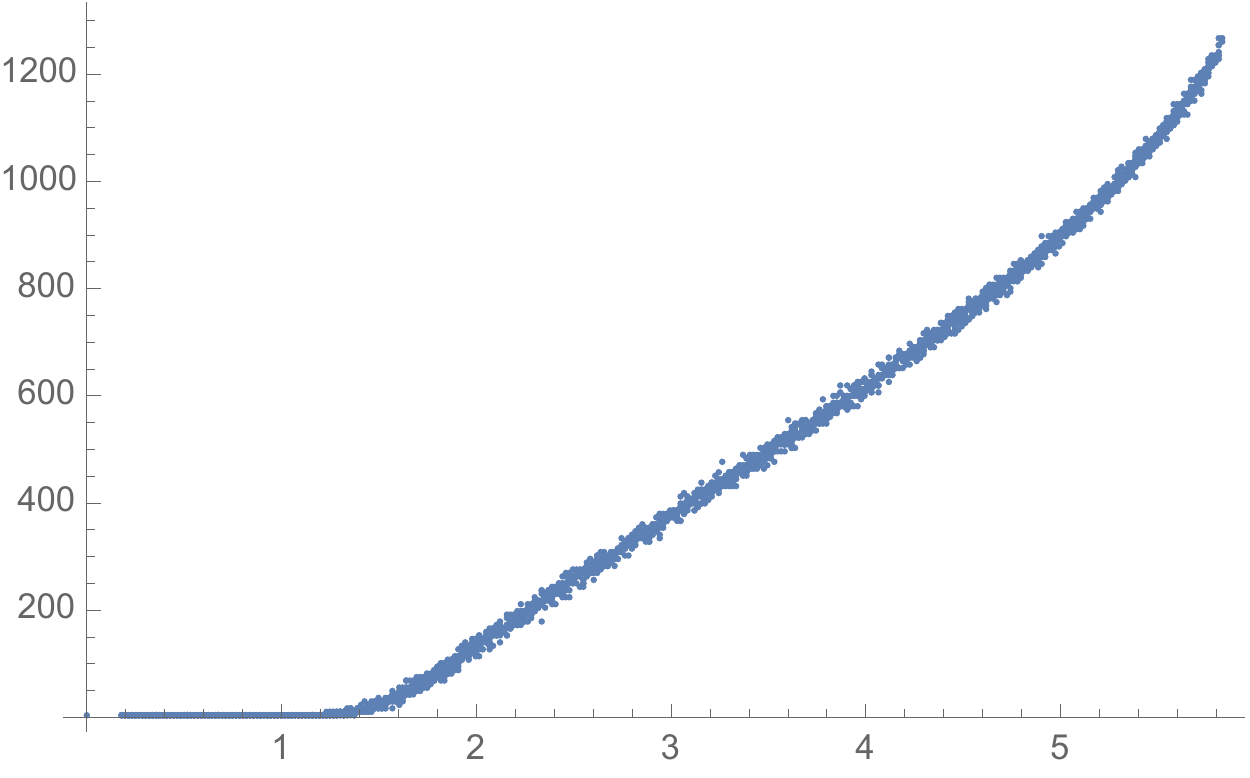}
\caption{\label{fig:cubicNodal}Number of nodal domains for random cubic graph}
\end{figure}

\begin{figure}
\centering
\includegraphics[width=0.7\textwidth]{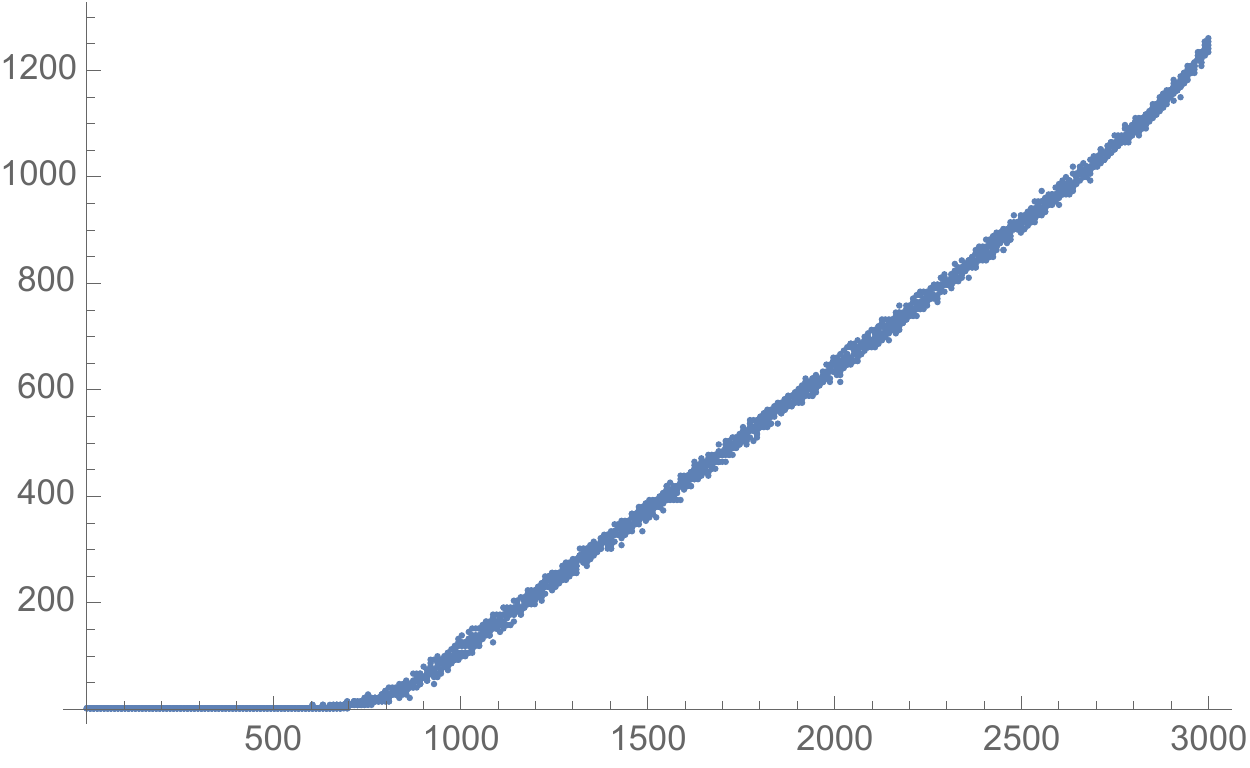}
\caption{\label{fig:cubicOrdinal}Number of nodal domains for random cubic graph by ordinal number of eigenvalue}
\end{figure}

\begin{figure}
\centering
\includegraphics[width=0.7\textwidth]{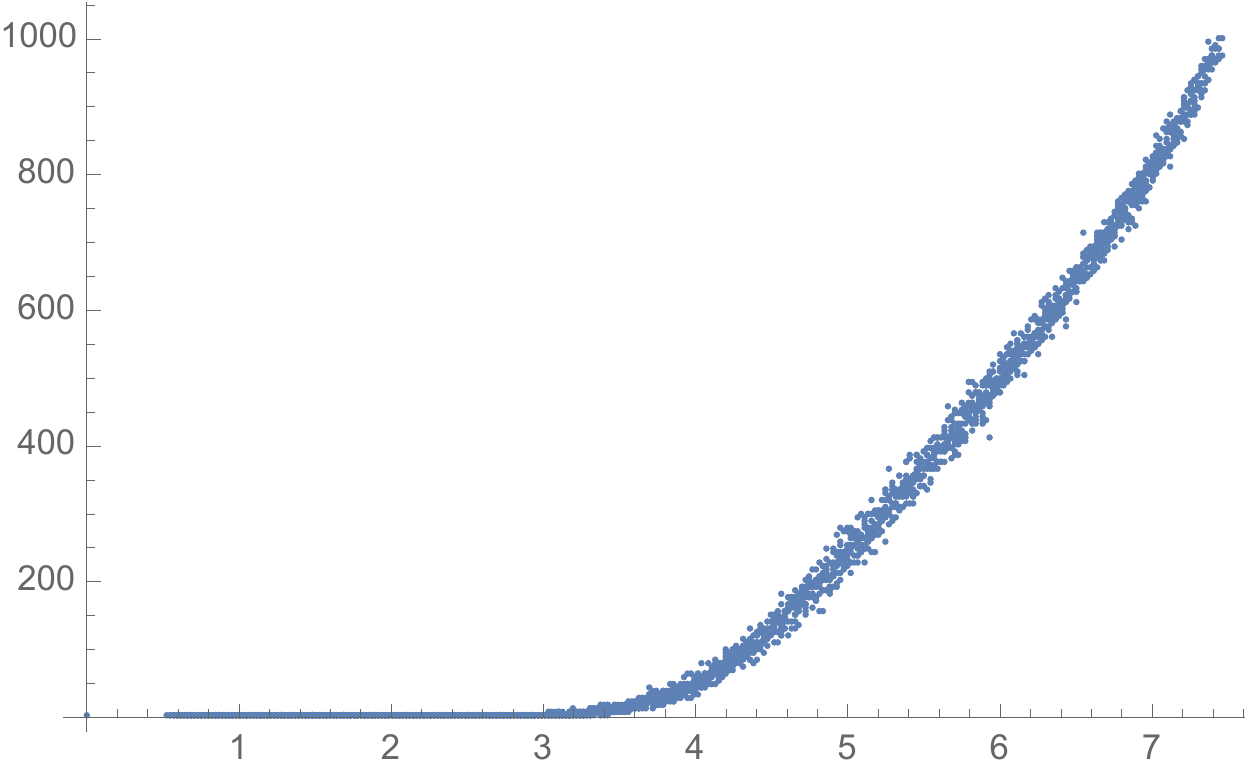}
\caption{\label{fig:quarticNodal}Number of nodal domains for random quartic graph}
\end{figure}

\begin{figure}
\centering
\includegraphics[width=0.7\textwidth]{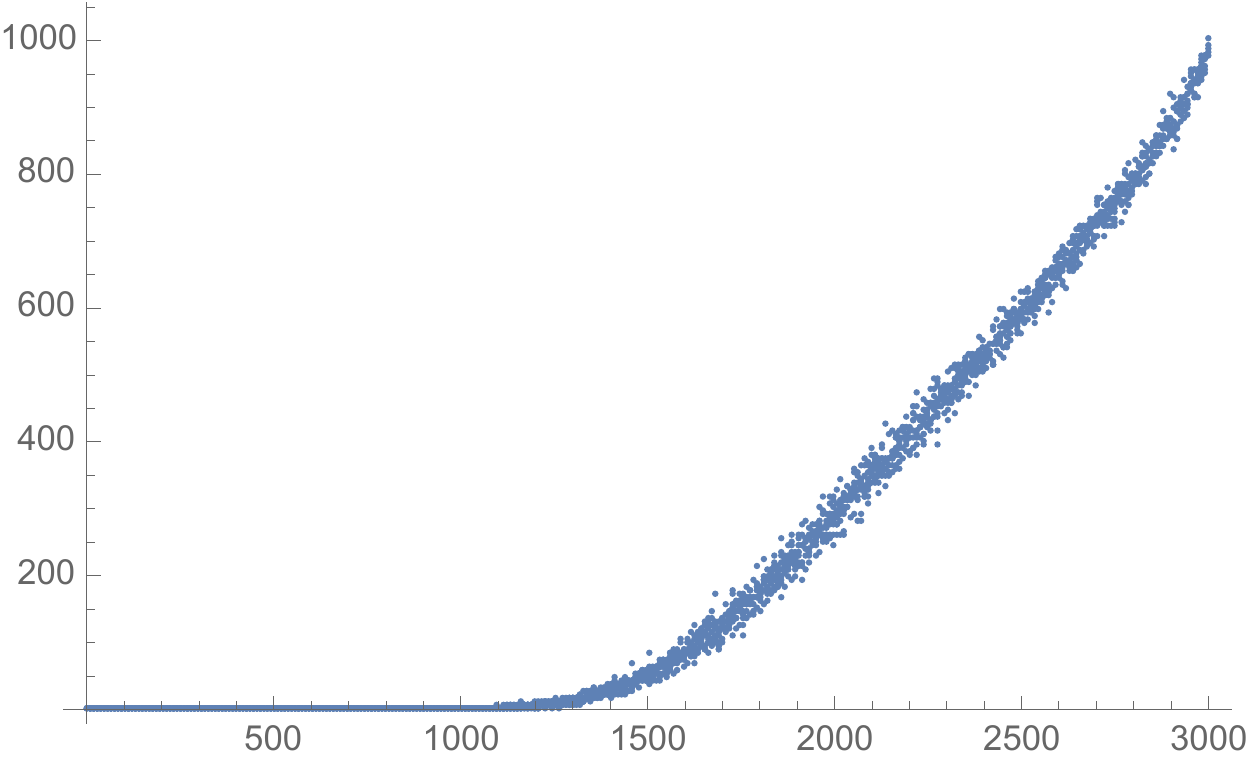}
\caption{\label{fig:quarticOrdinal}Number of nodal domains by ordinal number of eigenvalue for random quartic graph by ordinal number of eigenvalue}
\end{figure}

\begin{figure}
\centering
\includegraphics[width=0.7\textwidth]{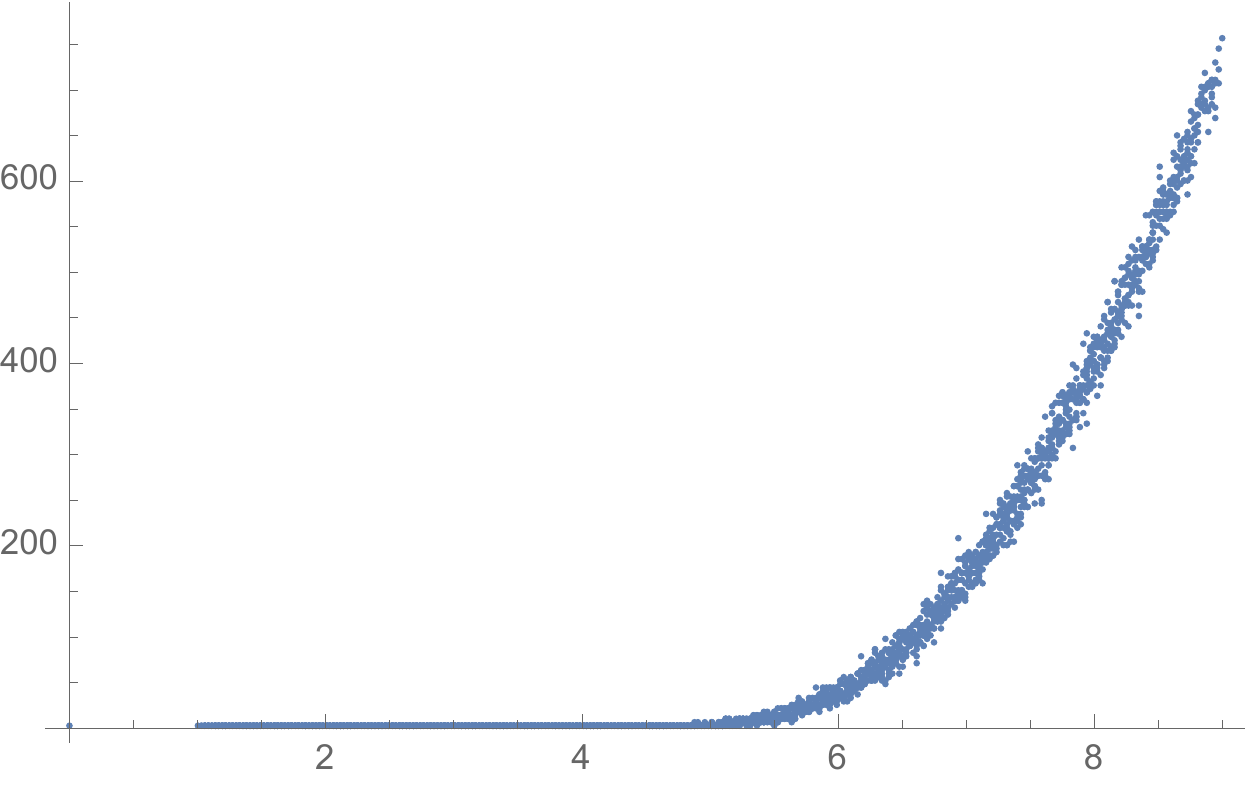}
\caption{\label{fig:quinticNodal}Number of nodal domains for random quintic graph}
\end{figure}

\begin{figure}
\centering
\includegraphics[width=0.7\textwidth]{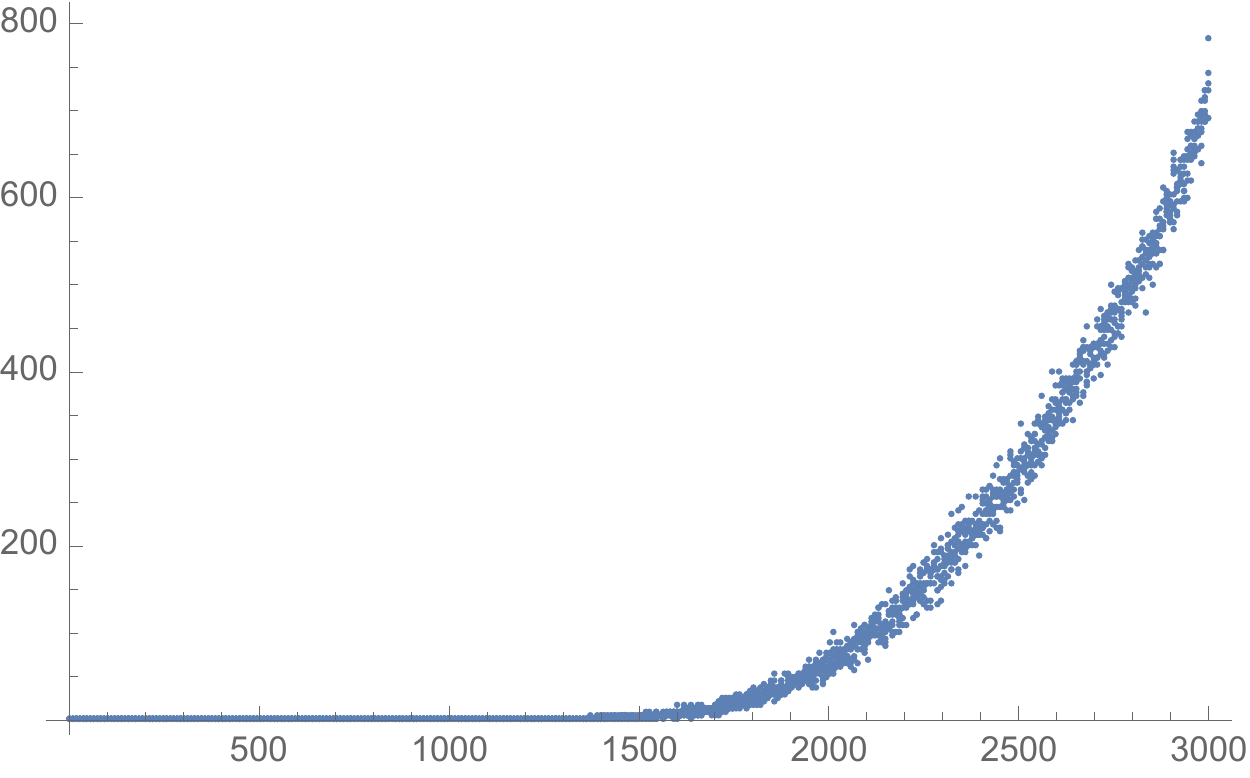}
\caption{\label{fig:quinticOrdinal}Number of nodal domains for random quintic graph by ordinal number of eigenvalue}
\end{figure}

\begin{figure}
\centering
\includegraphics[width=0.7\textwidth]{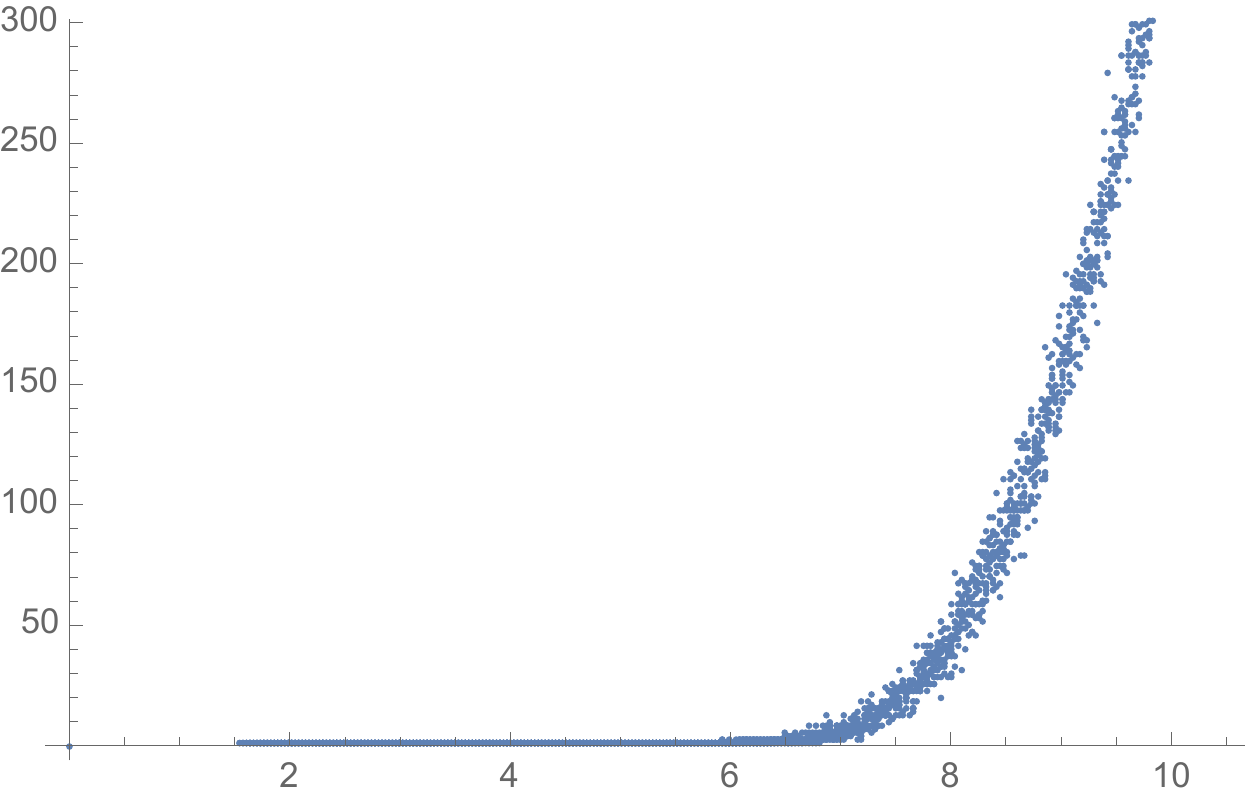}
\caption{\label{fig:sexticNodal} Number of nodal domains for random sextic graph}
\end{figure}

\begin{figure}
\centering
\includegraphics[width=0.7\textwidth]{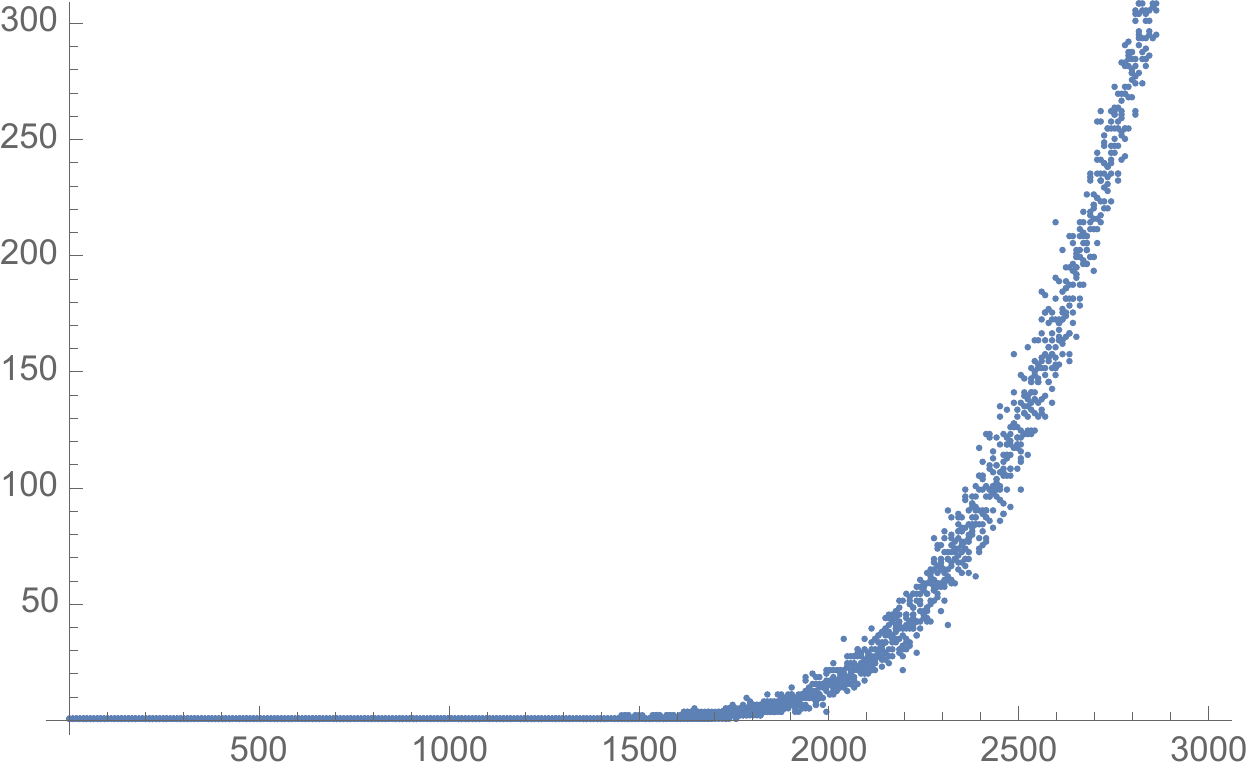}
\caption{\label{fig:sexticOrdinal} Number of nodal domains for random sextic graph by ordinal number of eigenvalue}
\end{figure}

In this case, it turns out that Theorem \ref{yauthm} is about as far from sharp as it can be. The graphs below, show the numbers of positive strong nodal domains for $3, 4, 5,$ and $6$-regular graphs. The number of vertices in our test graph equal $3000,$ and you see that, except in the cubic case, the number of strong nodal domains in this case is $O(1)$ (more precisely, it is equal to $2.$) See Figures \ref{fig:cubicNodal}, \ref{fig:quarticNodal} , \ref{fig:quinticNodal}, \ref{fig:sexticNodal}, and also Figures \ref{fig:cubicOrdinal}, \ref{fig:quarticOrdinal}, \ref{fig:quinticOrdinal}, \ref{fig:sexticOrdinal}. In the first batch of figures, the abscissa is labeled by the  numerical eigenvalue, in the second, by the ordinal number of the eigenvalue. We see that there is a phase transition: in the low energy regime, there are exactly two nodal domains (the positive and the negative), and then the number starts to grow, in a fairly fixed pattern (concave down, but close to linear). 
It is also quite clear that once the degree is above $4,$ the phase transition occurs at \textbf{exactly} the middle of the spectrum (corresponding to the $0$ eigenvalue of the adjacency matrix) - in fact, as shown in \cite{band2007nodal}, this transition cannot happen any later than that. Because of the low-energy flat piece, it seems clear that no general lower bound is possible. However, exactly such a lower bound is stated by Xu and Yau (Theorem 1.3 in \cite{xu2012nodal}:
\begin{theorem}[Xu-Yau, \cite{xu2012nodal}]
Let $f$ be an eigenfunction corresponding to $\lambda_k,$ which is zero on exactly $z$ vertices. Then the number of strong nodal domains of $f$ is at least $k + r - l - z,$ where $r$ is the multiplicity of $\lambda_k,$ and $l$ is the minimal number of edges that need to be removed from $G$ in order to turn it into a forest.
\end{theorem}
Now, experimentally, $z$ is always $0$ for our random graphs, and since the graphs are connected, $l = d n/2 - n+1,$ so in our case the Xu-Yau lower bound is at least $\max(1,k + n + 1 - d n/2) =\max(1, k + n(1-d/2) + 1).$ 

Interestingly, this lower bound seems is not bad for cubic graphs, but once $d \geq 4$ it becomes vacuous.

\subsection{Random planar maps}
\begin{figure}
\centering
\includegraphics[width=0.7\textwidth]{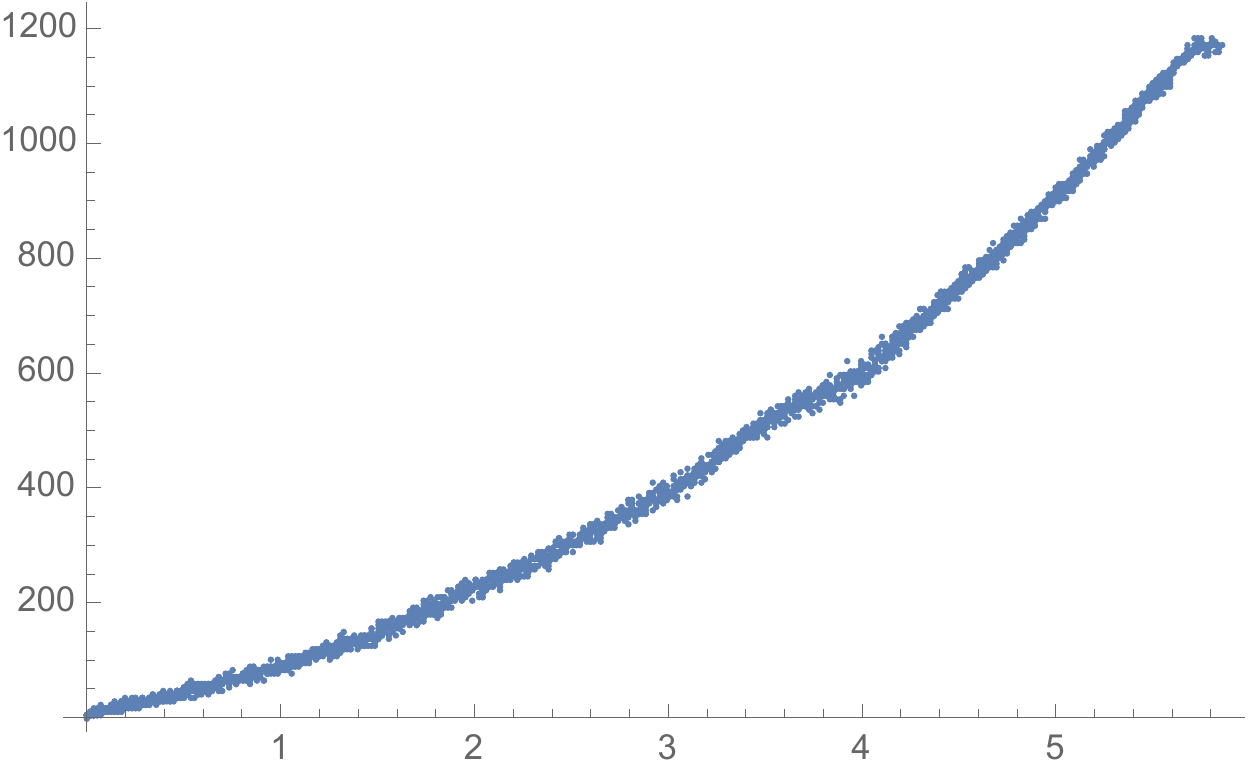}
\caption{\label{fig:planarNodal} Number of nodal domains for random $3$-conneceted cubic planar map}
\end{figure}

The statistic of nodal domains of random $3$-connected planar cubic maps (see Figure \ref{fig:planarNodal}) - the graph is (more or less) convex, and the maximum seems to be very close to where it is for random cubic graphs. What is quite visible, however, is that there is no longer the low-energy flat piece. Planar graphs were studied by Lin \emph{et al} in \cite{lin2013nodal}, and they show that in the planer $3$-connected case, the number of nodal domains of the $k$-th eigenvalue is at most $6k-34.$ Of course, this estimate is only of interest for $k < V(G)/6,$ and it is quite clear that in that range (for \emph{random} planar maps) the coefficient $6$ is a massive overestimate of the real growth rate.
\subsection{Random Voronoi triangulations}
\begin{figure}
\centering
\includegraphics[width=0.7\textwidth]{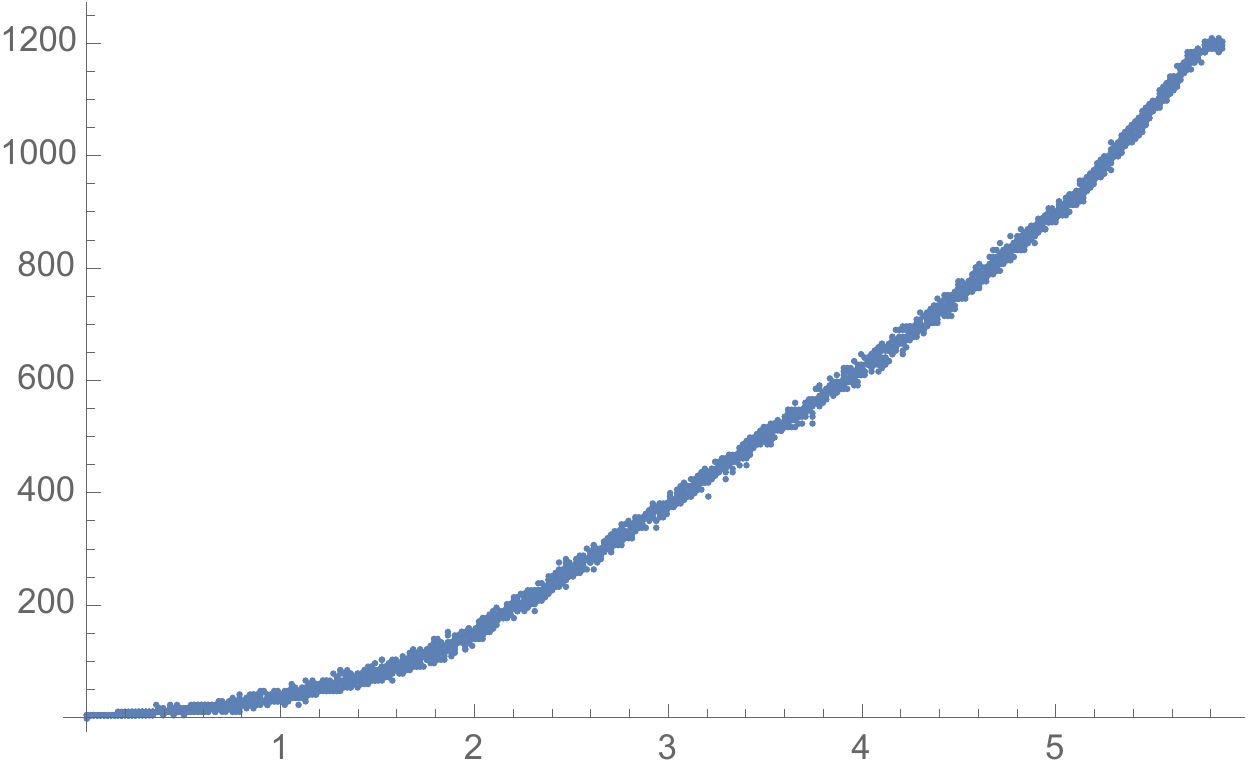}
\caption{\label{fig:voronoiNodal} Number of nodal domains for the Voronoi diagram of a random point set on the sphere}
\end{figure}

We show the graph in Figure \ref{fig:voronoiNodal}. The noteworthy differences between it and the random cubic planar maps seem to be a shallower low energy segment, and a much lower maximal number (in fact, it looks like the curve flattens out again at the end, a portent for what happens for Delaunay triangulations in the next figure.
\subsection{Random Delaunay Triangulations}
\begin{figure}
\centering
\includegraphics[width=0.7\textwidth]{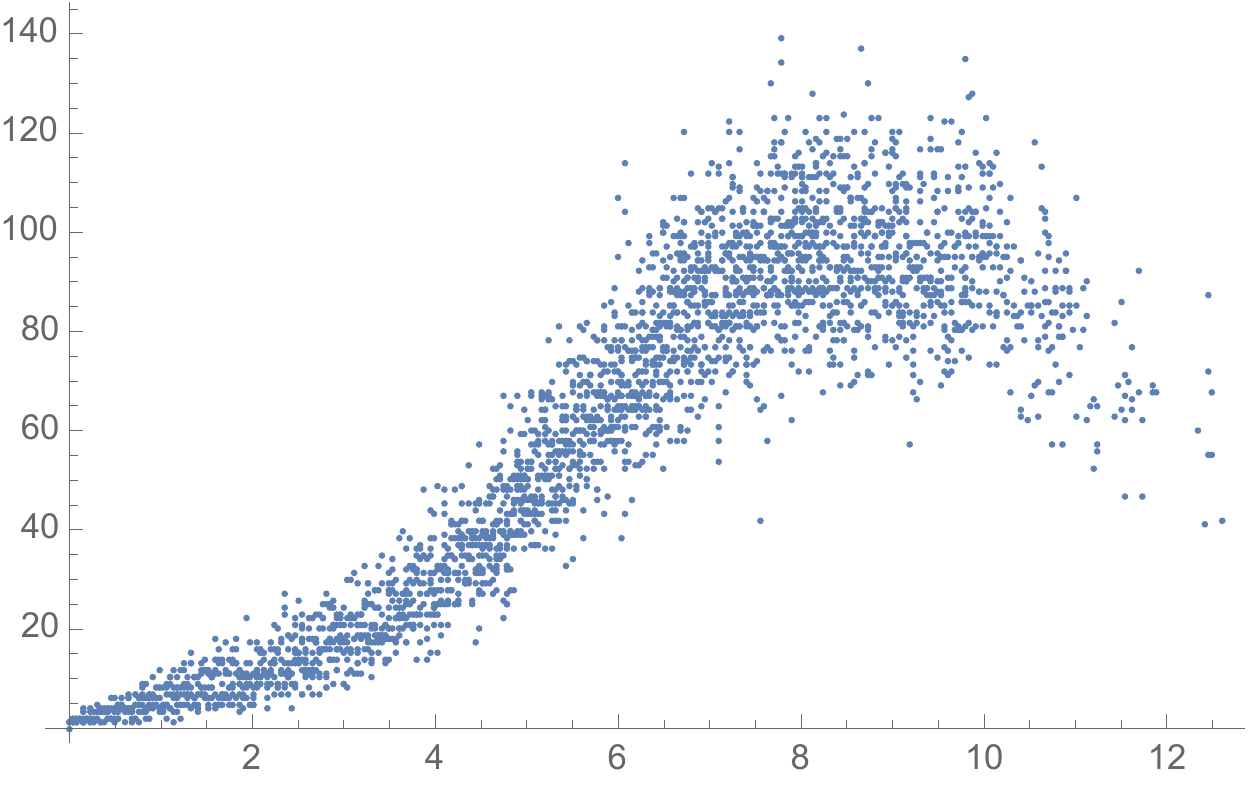}
\caption{\label{fig:delaunayNodal} Number of nodal domains for Delaunay triangulation of a random point set on the sphere}
\end{figure}

The strangest nodal domain distribution comes from random Delaunay triangulations of random point sets (see Figure \ref{fig:delaunayNodal}). It can be seen that the distribution is much less concentrated than in all other cases (presumably this is so because the graphs are not regular), but more interesting is the very noticeable \emph{drop} in the number of nodal domains at high energies.
\bibliographystyle{plain}
\bibliography{spectra}
\end{document}